\documentclass{article}
\usepackage[a4paper, total={6in, 8in}]{geometry}
\usepackage[utf8]{inputenc}
\usepackage[T1]{fontenc}
\usepackage{amsmath}
\usepackage{amsfonts}
\usepackage{amssymb}
\usepackage{amsthm}
\usepackage[pdfencoding=auto]{hyperref}
\usepackage{listings}
\usepackage{xcolor}
\usepackage[noadjust]{cite}
\usepackage{pgf,tikz,pgfplots}
\pgfplotsset{compat=1.15}
\pgfplotsset{ticks = none}
\usepackage{mathrsfs}
\usetikzlibrary{arrows}
\usepackage{verbatim}
\usepackage{cleveref}
\usepackage{enumitem}
\usepackage{yhmath}
\usepackage{mathtools}

\theoremstyle{plain}
\newtheorem{theorem}{Theorem}[section]
\newtheorem{lemma}[theorem]{Lemma}
\newtheorem{corollary}[theorem]{Corollary}
\newtheorem{proposition}[theorem]{Proposition}

\newtheorem{observation}[theorem]{Observation}

\newtheorem{maintheorem}{Theorem}

\newtheorem{maincorollary}[maintheorem]{Corollary}

\theoremstyle{definition}
\newtheorem{remark}[theorem]{Remark}
\newtheorem{definition}[theorem]{Definition}
\newtheorem{example}[theorem]{Example}

\newcommand*{\N}{\mathbb{N}}
\newcommand*{\Z}{\mathbb{Z}}

\newcommand*{\Cop}{\mathcal{C}_{\text{op}}^-}

\newcommand*{\lk}{\operatorname{link}}
\newcommand*{\hi}{\hat{\imath}}
\newcommand*{\hj}{\hat{\jmath}}

\newcommand*{\uak}{\mathcal{U}_A(k)}

\newcommand*{\CC}{\mathcal{CC}}
\newcommand*{\GAk}{\mathcal{G}_A(k)}
\newcommand*{\chev}{\operatorname{Chev}}
\newcommand*{\Ad}{\operatorname{Ad}}
\newcommand{\Ao}{\mathring A}

\newcommand*{\height}{\operatorname{ht}}
\newcommand*{\id}{\operatorname{Id}}
\newcommand*{\abar}{\overline{\alpha}}
\newcommand{\ophiJ}{\mathring \Phi_J^+}
\newcommand{\ophiK}{\mathring \Phi_K^+}
\newcommand{\stab}{\operatorname{Stab}}

\usepackage{todonotes}
\usepackage{parskip}
\newcommand{\footremember}[2]{%
    \footnote{#2}
    \newcounter{#1}
    \setcounter{#1}{\value{footnote}}%
}

\title{High-dimensional expanders from Kac--Moody--Steinberg groups}
\author{%
  Laura Grave de Peralta\footremember{alley}{Research Institute in Mathematics and Physics, UCLouvain, Chemin du Cyclotron 2, B--1348 Louvain--la--Neuve, Belgium}%
  \and Inga Valentiner-Branth\footremember{trailer}{Department of Mathematics: Algebra and Geometry, Ghent University, Krijgslaan 281 – S25, 9000 Gent, Belgium}%
  \footnote{Both authors are supported by the FWO and the F.R.S.--FNRS under the Excellence of Science (EOS) program (project ID~40007542).}
}
\date{January 31, 2025}

\begin{document}

\maketitle

\begin{abstract}
High-dimensional expanders are a generalization of the notion of expander graphs to simplicial complexes and give rise to a variety of applications in computer science and other fields. 
We provide a general tool to construct families of bounded degree high-dimensional spectral expanders. Inspired by the work of Kaufman and Oppenheim, we use coset complexes over quotients of Kac--Moody--Steinberg groups of rank $d+1$, $d$-spherical and purely $d$-spherical. We prove that infinite families of such quotients exist provided that the underlying field is of size at least 4 and the Kac--Moody--Steinberg group is 2-spherical, giving rise to new families of bounded degree high-dimensional expanders.
In the case the generalized Cartan matrix we consider is affine, we recover the construction in \cite{o2022high} (and thus also \cite{kaufman2018construction}) by considering Chevalley groups as quotients of affine Kac--Moody--Steinberg groups. Moreover, our construction applies to the case where the root system is of type $\tilde{G}_2$, a case that was not covered in earlier works. 
\end{abstract}\section{Introduction}
Expander graphs can be described as being highly connected yet sparse, despite how contradictory this may seem. The notion of expander graphs can be traced back to the 70s and proved useful in computer science, for example for constructing networks that are reliable as well as cost-effective. A lot of research and work has been devoted to expander graphs, involving not only mathematicians. We refer to \cite{lubotzky2012expander} for a survey.

The notion of expansion has been extended to higher dimension, yielding several, non-equivalent, definitions of expansion for simplicial complexes. Each has its own advantages and disadvantages. Among the most used generalizations are geometric expanders, topological expanders, coboundary expanders, cosystolic expanders and, finally, spectral expanders. These different notions are summarized in the survey \cite{lubotzky2018high}. High-dimensional expanders are connected to theoretical computer science, in particular to locally testable codes and to quantum LDPC codes. They are also related to matroid bases and the resolution of the Mihail-Vazirani conjecture (\cite{mihail1989expansion}). For more perspectives on high-dimensional expanders, see \cite{gotlib2023no}.

Similar to the study of expander graphs, one main goal is to construct infinite families of high-dimensional expanders, with growing size but with constant bounded degree. The first construction of such families was presented in 2005 by Lubotzky, Samuels and Vishne. The construction relies on the theory of Bruhat--Tits buildings and gives rise to so-called Ramanujan complexes (\cite{LubSamVishRamanujan}). The next construction by Kaufman and Oppenheim from 2018 (\cite{kaufman2018construction}) uses the more elementary tools of coset complexes over explicit groups to construct spectral high-dimensional expanders. This construction was generalized four years later by O'Donnell and Pratt (\cite{o2022high}) to all Chevalley groups except those of type $G_2$. Our construction will, in some sense, again be a generalization of the two aforementioned ones. The main novelty we introduce is that our construction is based on Kac--Moody--Steinberg groups (KMS groups) whose definition is recalled in \Cref{section: KMS group and root subgroups}. We give a general tool to produce high-dimensional expanders for any nice family of finite quotients of the KMS groups. In particular, our approach allows us to define infinite families of high-dimensional expanders associated to the root system $G_2$.

In this paper, we work with the notion of local spectral expansion for simplicial complexes. For a simplicial complex $X$, we will say it is a local spectral expander (or HDX for short) if all underlying graphs of all links in $X$ have good expansion properties. 
In particular, we rely on Oppenheim's \emph{trickling down theorem} (see \cite{oppenheim2018local}). This result provides an easy criterion to determine expansion for a simplicial complex: given good connectivity properties of the simplicial complex $X$, we only need to check the expansion properties of the links in $X$ that have dimension 1, that is, those links which are graphs. 

This tool by Oppenheim has been used in \cite{kaufman2018construction}, where the authors prove that the coset complex associated to some elementary matrix groups over a ring $R$ and some well-chosen subgroups have local spectral expansion. A coset complex $\CC(G;(H_i)_{i=0}^d)$ over a group $G$ and subgroups $H_0,\dots,H_d \leq G$ is a pure, $d$-dimensional simplicial complex with vertex set $\bigsqcup_{i=0}^d G/H_i$ and the maximal faces are of the form $\{gH_0,\dots,gH_d\}$ for $g\in G$ (see \Cref{def: coset complex}). In \cite{o2022high}, O'Donnell and Pratt generalize the construction described by Kaufman and Oppenheim to construct local spectral high-dimensional expanders from Chevalley groups, also using the trickling down theorem.

We consider a further generalization of Chevalley groups $\chev_{A}(k)$. Chevalley groups are constructed from finite-dimensional Lie algebras encoded by a Cartan matrix $A$ and a field $k$. A similar construction for \emph{generalized} Cartan matrices yields groups associated to infinite-dimensional generalizations of the finite-dimensional simple Lie algebras called Kac--Moody algebras. More precisely, we look at Kac--Moody--Steinberg groups. They can be defined as fundamental groups of certain complexes of finite $p$-groups, where $p$ is a fixed prime number. The finite $p$-groups involved in the construction are the positive unipotent subgroups of basic Levi subgroups of spherical type in a 2-spherical Kac--Moody group over a finite field. A specific Kac--Moody--Steinberg group first appeared in \cite{ershov2010property} and \cite{ErshovG-SGrps} as an example of a Golod–Shafarevich group with property (T). The rank 3 case has been further investigated in \cite{PEC2022Hyperbolic}.

Kac--Moody groups and thus also Kac--Moody--Steinberg groups are determined by the underlying root system which can be described using a generalized Cartan matrix $A = (A_{i,j})_{i,j \in I}$ or a Dynkin diagram on $|I|$ nodes, which is a labelled, oriented graph where the edges are determined by the entries of $A$ (see \Cref{section: KMS group and root subgroups}). We call $|I|$ the rank of $A$. A subset $J\subseteq I$ is called spherical if the submatrix $(A_{i,j})_{i,j \in J}$ of the generalized Cartan matrix (GCM) is of finite type, which in particular implies that the associated root system is finite. A GCM is called $n$-spherical if every subset of size $n$ is spherical.

A 2-spherical Kac--Moody--Steinberg group $\mathcal{U}_A(k)$ over a field $k$ (or KMS group for short) is generated by its root subgroups $U_i$ (also called \emph{local groups}), associated to the simple root $\alpha_i$ for $ i \in I$. We assume for our construction that the GCM $A$ is of rank $d+1$, $d$-spherical and purely $d$-spherical (i.e. every spherical subset is contained in a spherical subset of size $d$). This includes the affine types. We set $U_{\hi} = \langle U_j \mid j \in I \setminus \{ i \} \rangle$.

Our construction exploits the fact that, in general, Kac--Moody--Steinberg groups have infinite families of finite quotients in which the local groups embed. From such data, we construct a coset complex and prove it has high-dimensional spectral expansion. 

Our main result is the following.

\begin{maintheorem}[\Cref{thm: main thm}] \label{thm: A}

Let $G$ be a finite group, $\phi$ a surjective morphism from a rank $d+1$, $d$-spherical, purely $d$-spherical Kac--Moody--Steinberg group $\mathcal{U}_A(k)$ to $G$ and let $H_i$ be the images of the $U_{\hi}$ in $G$. We consider the coset complex $X$ of $G$ with respect to the subgroups $\{H_i\}_{i \in I}$:
    $$X = \CC(G; (H_i)_{i\in I}).$$
If the map $\phi: \uak \to G$ satisfies
\begin{enumerate}[label=(\alph*)]

    \item $\phi$ is injective on the local groups $U_J := \langle U_j \mid j \in J \rangle$ for each spherical subset $J \subseteq I$
    \item for every spherical subsets $J,K \subseteq I$ we have
    $$\phi(U_J) \cap \phi(U_K) = \phi(U_J \cap U_K),$$
\end{enumerate}

then $X$ is a pure $d$-dimensional simplicial complex. Moreover, there exists constants $ \gamma >0$, where $\gamma $ only depends on the field $k$ and the GCM $A$ (see \Cref{thm: main thm} for the exact values), and $\Delta >0$ such that the following hold.
\begin{enumerate}
    \item For any face $\sigma$ of dimension $i \leq d-2$, the link of $\sigma$, $\lk_X(\sigma)$, is connected. In particular, $X$ is connected.
    
    \item For all $v \in X(0)$: $|\{\sigma \in X(d): v \subseteq \sigma \}| \leq \Delta$.
    
    \item\label{main thm: local expansion} For any $\sigma \in X(d-2): \lambda_2(\lk_X(\sigma)) \leq \gamma$, where $\lambda_2(\lk_X(\sigma))$ is the second-largest eigenvalue of the weighted random walk operator on $\lk_X(\sigma)$.
    
    \item $G$ acts sharply transitively on the maximal faces of $X$.
    \end{enumerate}
In particular, if $k$ is chosen such that $\gamma \leq \frac{1}{d}$ then $X$ is a $\left( \frac{\gamma}{1-(d-1)\gamma} \right)$ -- spectral high dimensional expander.

\end{maintheorem}

In case the GCM is 3-spherical and the finite field $k$ is of size $\geq 5$, it is known that $\uak$ is residually $p$ (where $p$ is the characteristic of the field $k$). Thus, we get infinite families of finite quotients such that the projections are injective enough to satisfy the conditions in \Cref{thm: A}. In the case the GCM is 2-spherical, even though the group is not residually $p$, we can also construct families of finite quotients satisfying the requirements of \Cref{thm: A} (this is done in \Cref{subsection: quotients of KMS groups}). 

Another way to get families of finite quotients of certain KMS groups is by starting with a spherical root system $\Ao$ and its affinization $A$. Then we get a natural morphism from the KMS group $\uak$ to the Chevalley group associated to $\Ao$ over the polynomial ring in one variable $\chev_{\Ao}(k[t])$ which is injective on the local groups. Fixing a family of irreducible polynomials of growing degree $f_m$ and passing to the quotient $\chev_{\Ao}\left(k[t]/(f_m)\right)$ gives rise to the desired family of maps. 
The latter choice of quotients yields high-dimensional expanders that are very similar to the ones considered in \cite{kaufman2018construction} and \cite{o2022high}. However, there are some differences with the techniques we use for bounding the degree in the complexes.

To summarize, we get the following result. 

\begin{maincorollary}[\Cref{corol: main corollary}]
     Let $A$ be a 2-spherical GCM of rank at least 3 and let $k$ be a finite field of size at least 4. 
    There exists an infinite family of finite index subgroups $N_1, N_2,\dots \unlhd \uak$, with $[\uak: N_i] \to +\infty$ such that if we set $G_i:= \uak / N_i$, $\phi_i: \uak \to G_i$ the quotient map, $H_{j}^{(i)}:= \phi_i(U_{\hj}), j \in I$ the images of the local groups and $X_i:= \CC(G; (H_j^{(i)})_{j\in I})$, then $\left( X_i \right)_{i \in \N}$ is a family of bounded degree high-dimensional spectral expanders.
\end{maincorollary}

This paper is organized as follows. First, in \Cref{section: preliminaries}, we introduce the necessary combinatorial tools needed for our result. We state the notion of high-dimensional spectral expansion in \Cref{subsection: HDX definitions}, then we introduce coset complexes in \Cref{subsection: coset complexes}, and finally we talk about chamber systems and their relation to coset complexes in \Cref{subsection: chamber systems}. \Cref{section: KMS group and root subgroups} introduces the necessary Kac--Moody theory. We discuss Kac--Moody--Steinberg groups and their root groups in \Cref{subsec: KMS groups}. In \Cref{subsec: action on buildings} we take a look at the action of Kac--Moody groups on (twin) buildings and its consequences. We combine the results to show the existence of finite quotients of Kac--Moody--Steinberg groups in \Cref{subsection: quotients of KMS groups}. This is necessary for the construction of families of expanders. \Cref{section: main theorem} contains our main result, that is  proved in \Cref{subsection: connectivity links,subsection: bounded degree ,subsection: link expansion,subsection: symmetries}. In \Cref{section: Chevalley } we discuss how our construction relates to constructions using Chevalley groups. \Cref{subsection: quot in Chev: Def and Prelim} provides the ideas and basic steps of the construction. The results are proved  in \Cref{subsection: quot in Chev: the general case}. In \Cref{subsection: quot in Chev: sl3}, we give an explicit example of our construction in the case where the Chevalley group is $\operatorname{SL}_3$. Finally, we make general remarks and comments about our construction in \Cref{section : context}. In \Cref{subsection: quot in chev: comparison}, we compare our construction to previous results by O'Donnell and Pratt \cite{o2022high} and Kaufmann and Oppenheim \cite{kaufman2018construction}. In \Cref{subsection : generalization}, we restate our main result in a more general setting that does not rely on Kac--Moody theory.

\section*{Acknowledgements}
We are grateful to Pierre-Emmanuel Caprace, Tom De Medts and Timothée Marquis for introducing us to this question and topic and their comments. We especially want to thank Pierre-Emmanuel Caprace for sharing his ideas about how to combine Kac--Moody--Steinberg groups and high-dimensional expanders. We are grateful to Jari Desmet, Alex Loué, François Thilmany, Ferdinand Ihringer and Corentin Le Coz for fruitful discussions.

We thank both referees for their suggestions that substantially improved the exposition of the paper. We further want to thank one of the referees for their suggestion to use a different approach to prove the link expansion. The new approach improved the bounds compared to the ones we achieved with our original method.

Both authors are supported by the FWO and the F.R.S.--FNRS under the Excellence of Science (EOS) program (project ID~40007542).

\section{Preliminaries}\label{section: preliminaries}
We will write $\N$ for the set of non-negative integers and $\N^*$ for the set of positive integers.

\subsection{High-dimensional spectral expanders (HDX)}\label{subsection: HDX definitions}

Expander graphs are sparse yet highly connected graphs. They can be defined in several ways (see \cite{lubotzky1994discrete} for an overview), but all definitions are equivalent. In the case of simplicial complexes of dimension $d\geq 2$, this is not the case anymore and several notions exist: coboundary expanders, cosystolic expanders, topological expanders, geometric expanders, spectral expanders etc. (see \cite{lubotzky2018high}).  In this paper, we will consider the notion of $\lambda$\emph{-spectral expansion} for families of simplicial complexes.

A \emph{simplicial complex} $X$ with vertex set $X(0)$ is a collection of subsets of $X(0)$ satisfying:

\begin{itemize}
    \item $\{v\} \in X$ for every $v \in X(0)$;
    \item If $\sigma \in X$ and $\tau \subseteq \sigma$ then $\tau \in X$.
\end{itemize}

For $i=-1,0,1,2\dots$ we denote by $X(i)$ the set of subsets of $X(0)$ of size $i+1$, where we follow the convention of identifying a vertex $v$ with the set containing it $\{v\}$. An element $\sigma \in X(i)$ is called a \emph{face of dimension $i$} or a \emph{simplex}, and we will usually refer to elements in $X(0)$ as \emph{vertices} and elements in $X(1)$ as \emph{edges} of the simplicial complex. If there exists $d \in \N$ such that $X(k)=\varnothing$ for all $k > d$ and $X(d) \neq \varnothing$, then $X$ is said to be finite dimensional and the \emph{dimension} of $X$ is $d$. A simplicial complex $X$ of dimension $d$ is \emph{pure} if every face is contained in a maximal face of dimension $d$. For a simplex $\sigma$ in a $d$-dimensional simplicial complex, we define $\deg(\sigma) = \lvert\left\{\tau \in X(d) \mid \sigma \subseteq \tau \right\} \rvert$ and we say that $X$ is of $\Delta$-bounded degree if every vertex  belongs to at most $\Delta$ maximal faces, i.e. $\deg(\sigma) \leq \Delta$ for all non-empty $\sigma \in X$. 
We will only work with finite, pure simplicial complexes of bounded degree.

For a face $\sigma \in X$, we define its \emph{link} in $X$, denoted $\lk_X(\sigma)$, to be the set of simplices $\tau$ such that $\tau \cup \sigma \in X$ and $\tau \cap \sigma= \varnothing$, or equivalently, $\lk_X(\sigma)=\left\{ \tau \setminus \sigma \ | \ \tau \in X, \ \sigma \subseteq \tau  \right\}$. The \emph{$k$-skeleton} of a simplicial complex $X$ is $X(0) \cup \dots \cup X(k)$, the collection of all faces of dimension up to $k$. 

We are particularly interested in the 1-skeleton of the complex $X$ (that is, the underlying graph) and the 1-skeleton of the links, $\lk_X(\sigma)$ for $\sigma \in X(i)$ and $i \leq d-2$. We say that a simplicial complex is \emph{connected} if its 1-skeleton is a connected graph. We denote by $K_{\sigma}$ the weighted  1-skeleton of $\lk_X(\sigma)$ (omitting $X$ and the weight considered as it will be clear from the context)  where a \emph{weight} on the simplicial complex $X$ is defined as follows. 

\begin{definition}[{\cite[Definition 2.1]{oppenheim2018local}}]
    Let $X$ be a $d$-dimensional simplicial complex. A \emph{weight} $w$ on $X$ is a function $w: \bigcup_{-1 \leq k \leq d} X(k) \to \mathbb{R}^+$. The weight function is called \emph{balanced} if for every $-1 \leq k \leq d$ and every $\tau \in X(k)$, we have  
    \[ 
    \sum_{\sigma \in X(k+1), \tau \subset \sigma} w(\sigma) = w(\tau).
    \]

    A simplicial complex $X$ with a balanced weight function will be called a weighted simplicial complex.
\end{definition}

We will work with a particular choice of weight. The idea here is to assign weight 1 to each maximal dimensional simplex and then extend the weight function such that it is balanced. We get the following weight function. For $\tau \in X{(k)}$, we define it as 
$$
w(\tau)=(d-k) !\left|\left\{\sigma \in X{(d)}: \tau \subseteq \sigma\right\}\right| .
$$

For every $-1 \leq k \leq d-2$, $\tau \in X(k)$, the one-skeleton of $\lk_X(\tau)$ will inherit the weight from $X$. Explicitly, given $\tau \in X(k)$ and $\{u, v\} \in \lk_X(\tau)(1)$, we have
$$
w_{\tau}(\{u, v\})=(d-k-2) !\left|\left\{\sigma \in \lk_X(\tau)(d-k-1):\{u, v\} \subseteq \sigma\right\}\right|.
$$
\begin{comment} (Old version)
We will write $\lambda_2(K_{\sigma})$ for the second-largest eigenvalue of the weighted random walk operator $M$ on the weighted graph $K_{\sigma}$ seen as an operator on $\ell^2(X(0))$. That is, we consider the random walk with probability of choosing an edge proportional to the weight of that edge. If $w_\sigma(\{u,v\})$ is the weight of the edge $\{u,v\}$ in $\lk(\sigma)$, then $$M_{u,v} =  \frac{w_\sigma(\{u,v\})}{\sum_{w \sim u} w_\sigma(\{u,w\})}$$ and this defines the random walk on the links.
\end{comment}

Next, we consider the weighted random walk on $K_\sigma$, that is we perform a random walk where the probability of choosing an edge is proportional to the weight of that edge. If $w_\sigma(\{u,v\})$ is the weight of the edge $\{u,v\}$ in $\lk(\sigma)$, then $$M_{u,v} =  \frac{w_\sigma(\{u,v\})}{\sum_{w \sim u} w_\sigma(\{u,w\})}$$ defines the weighted random walk operator $M$ acting on $\ell^2(\lk_X(\sigma)(0))$. We will write $\lambda_2(K_\sigma)$ for the second largest eigenvalue of $M$. Recall that $M$ has largest eigenvalue equal to 1.

If we consider the link of a simplex $\tau$ of dimension $d-2$, the random walk with respect to this choice of weight, for $X$ and thus for $\lk_X(\tau)$, is simply the standard random walk (as seen on an unweighted graph).

The definition of $\lambda$-spectral expansion we consider is the following, where $\lambda$ is a positive real number.

\begin{definition}(\cite{oppenheim2018local,kaufman2018construction})\label{def:HDX}
    A $d$-dimensional, pure, finite simplicial complex $X$ is said to be a $\lambda$-spectral (one-sided) high-dimensional expander if $\lambda_2(K_{\sigma}) \leq \lambda$ for all faces $\sigma$ of $X$ of dimension at most $d-2$. 
\end{definition}

We will write $\lambda$-HDX or simply HDX throughout the paper for shorter notation.

The \emph{spectral} notion of expansion considered here is a natural one to consider for simplicial complexes as it implies optimal mixing of high order random walks (see \cite{kaufman2017high, kaufman2018high, kaufman2023high}). Such a property is desirable since, if combined with good symmetries of the HDX, allows developing Low-Density-Parity-Check (LDPC) codes, among other applications (see e.g. \cite{EvraLDPC}).

By analogy with the graph theoretic setting, our goal is to construct families of $\lambda$-HDXs: that is a collection $\{X_n\}_{n \in \mathbf{N}}$ of $d$-dimensional pure simplicial complexes of bounded degree having an increasing number of vertices and such that they are $\lambda$-HDXs for $n$ large enough.

To prove a simplicial complex is a HDX in the sense of \Cref{def:HDX}, we will use the following result.

\begin{theorem}[\cite{oppenheim2018local}]\label{thm: criterion for HDX}
Let $X$ be a $d$-dimensional connected, pure simplicial complex and $0 < \gamma \leq \frac{1}{d}$ such that
\begin{enumerate}[label=(\alph*)]
    \item $X$ is connected;
    \item $\lk_X(\sigma)$ is connected for all $\sigma \in X(i)$ for $i \leq d-2$;
    \item $\lambda_2(K_\sigma) \leq \gamma $ for all $\sigma \in X(d-2)$.
\end{enumerate}
Then $X$ is a $\left( \frac{\gamma}{1-(d-1)\gamma} \right) $-- HDX.
\end{theorem}

This "trickling down" theorem states that, given good connectivity of the complex $X$, it suffices to prove expansion for the links of dimension $d-2$ (which are graphs).

\subsection{Coset complexes}\label{subsection: coset complexes}
In this section, we recall the definition and some known facts about coset complexes, following the summary on the same topic in \cite{o2022high}.

\begin{definition} \label{def: coset complex}
    Let $G$ be a group and $\mathcal{H}= (H_0,\dots,H_d)$ be a family of subgroups of $G$. The \emph{coset complex} $\CC(G; \mathcal{H})$ is the simplicial complex with 
    \begin{itemize}
        \item vertex set $\bigsqcup_{i=0}^d G/H_i$
        \item a set of vertices $\left\{ g_1H_{i_1},\dots,g_kH_{i_k} \right\}$ forms a $(k-1)$-simplex in $\CC(G; \mathcal{H})$ if and only if $\bigcap_{j=1}^k g_j H_{i_j} \neq \varnothing$.
    \end{itemize}
\end{definition}
The simplicial complex $\CC(G; \mathcal{H})$ is pure and $d$-dimensional. 

An equivalent characterization of a simplex in $\CC(G;\mathcal{H})$ is that a set of vertices $\left\{ g_1H_{i_1},\dots,g_kH_{i_k} \right\}$ forms a simplex if and only there exists $g \in G$ such that $gH_{i_j} = g_jH_{i_j}$ for all $j=1,..,k$. In particular, the maximal simplices of $\CC(G;\mathcal{H})$ are of the form 
$$\{gH_0,\dots, gH_d \}$$
for $g \in G$. \\

Note that the set of vertices is partitioned into $d+1$ subsets $G/H_i$ for $i = 0,\dots,d$ such that for any simplex $\sigma \in \CC(G;\mathcal{H})$ we have $|\sigma \cap G/H_i| \leq 1$ for all $i$. In that case, we say that $\CC(G;\mathcal{H})$ is a $(d+1)$-partite simplicial complex.

Examples of coset complexes are Coxeter complexes and certain buildings, for example when considering an isotropic semisimple algebraic group and its parabolic subgroups. 

\begin{definition} \label{def: type of a face}
    The \emph{type} of a simplex $\left\{ g_1H_{i_1},\dots,g_kH_{i_k} \right\}$ in $\CC(G;\mathcal{H})$ is the set of indices $\{i_1,\dots,i_k\} \subseteq \{0,\dots,d\}$.
\end{definition}

Note that the type of a simplex can be defined for any partite simplical complex. We observe the following.

\begin{lemma}\label{lem: action of G on CC}
    The group $G$ acts on $\CC(G;\mathcal{H})$ by left multiplication. This action is type-preserving (it does not change the type of a simplex) and transitive on maximal faces. The action is sharply transitive if $\bigcap_{i=0}^d H_i = \{1\}$.
\end{lemma}
\begin{remark} \label{rmk: converse of lemma}
    The converse of the above lemma is also true in the following sense. Given a pure, $d$-dimensional $d+1$-partite complex $X$ (i.e. a simplical complex with a type function on the vertices such that no two vertices of the same type are adjacent) and a group $G$ acting by type preserving automorphisms on $X$ such that the action on $X(d)$ is transitive, then $X$ isomorphic to a coset complex. If we fix one maximal simplex $\{v_0,\dots, v_d \} \in X(d)$ and denote $H_i = \operatorname{Stab}_G(v_i)$, then $X \cong \CC(G;(H_i)_{i=0}^d)$. A proof of this can be found for example in \cite[Proposition 5.5]{kaufman2023high}.
\end{remark}

\begin{definition} 
    Given a collection of subsets $\mathcal{H}=(H_0, \dots , H_d)$ of a group $G$ and a type $\varnothing \neq T \subseteq \{0,\dots,d\}$ we write 
    $$H_T:= \bigcap_{i \in T} H_i$$
    and set
    $$H_{\varnothing}:= \langle H_0,\dots,H_d \rangle \leq G.$$
\end{definition}

The following two well-known facts will be very useful. For detailed proofs, see \cite[Lemma 3.3 and Lemma 3.4]{harsha2019note}. 

\begin{proposition} \label{prop: connectivity of CC}
$\mathcal{C C}(G , \mathcal{H})$ is connected if and only if $H_{\varnothing}=G$.
\end{proposition}

\begin{proposition} \label{prop: links in CC} 
Let $\sigma$ be a face in $\CC (G , \mathcal{H})$ of type $T \neq \varnothing$. Then the link of $\sigma$ is isomorphic to the coset complex $\CC \left(H_T ,\left(H_{T \cup\{i\}}: i \notin T\right)\right)$.
\end{proposition}

\subsection{Chamber systems} \label{subsection: chamber systems}
In this section, we introduce the notion of chamber system, which is another approach to describe partite complexes like the coset complexes and buildings. For the exposition we follow \cite{ScharlauBuildings}. We introduce the notion of a chamber system being $m$-simply connected which is used later in \Cref{prop: U_A(k)=U+} to show properties of the KMS-groups. 
\begin{definition}{\cite[Definition 1.2.1]{ScharlauBuildings}}
    A chamber system over some index set $I$ is a set $\mathcal{C}$ together with equivalence relations $\sim_i$ on $\mathcal{C}$ for each $i \in I$. The elements of $\mathcal{C}$ are called chambers. 
\end{definition}
Similar to coset complexes, we can construct chamber systems using a group and a family of subgroups. 
\begin{proposition}\cite[Proposition 1.4.5]{ScharlauBuildings}
Let $G$ be a group, $B\leq G$ a subgroup and $G_i, i\in I$ a family of subgroups of $G$ containing $B$. Let $\mathcal{C}(G,B,G_i,i\in I) :=G/B$ together with the family of relations
\begin{align*}
    xB \sim_i yB \iff xG_i = yG_i.
\end{align*}
Then $\mathcal{C}(G,B,G_i,i\in I)$ is a chamber system over $I$.
\end{proposition}

\begin{definition}
\begin{enumerate}[label=(\roman*)]
    \item Let $(\mathcal{C}, \sim_i, i \in I)$ be a chamber system. A gallery in $\mathcal{C}$ is a chain $C_0 \sim_{i_1} C_1 \sim_{i_2} \dots \sim_{i_n} C_n$ such that $C_k \in \mathcal{C}, i_k \in I$ for all $k$. We call $(i_1,\dots, i_n)$ the type of the gallery.
    \item Two galleries $C=C_0 \sim_{i_1} \dots \sim_{i_n} C_n,  D= D_0 \sim_{j_1} \dots \sim_{j_k} D_k$ are called elementary $m$-homotopic if there exists a gallery $E = E_0 \sim_{h_1}\dots \sim_{h_\ell} E_\ell$, $0\leq k_1<k_2\leq n$ and a set $J \subseteq I$ with $\lvert J \rvert \leq m$ such that $\{h_0,\dots,h_\ell \} \subset J, \{i_{k_1+1},i_{k_1+2},\dots i_{k_2}\} \subseteq J$ and such that we have
    \begin{alignat*}{2}
        D = C_0 \sim_{i_1} \dots \sim_{i_{k_1}} \underset{\parallel}{C_{k_1}} & \qquad && \underset{\parallel}{C_{k_2}}  \sim_{i_{k_2+1}} \dots \sim_{i_n} C_n. \\
         E_0 & \sim_{h_1} \dots \sim_{h_\ell}  && E_\ell 
  \end{alignat*}
    \item We call two galleries $m$-homotopic if they can be connected by a finite sequence of elementary $m$-homotopies. 
    \item A chamber system is called $m$-simply connected if it is connected (i.e. any two chambers can be connected by a gallery) and every closed gallery is $m$-homotopic to a gallery of length 0.
    \end{enumerate}
\end{definition}
\begin{remark}{\cite[Chapter 1.3]{ScharlauBuildings}} \label{rmk: coco vs ch sys}
    Let $G$ be a group, $H_i, i\in I$ a family of subgroups over a finite index set $I$. Set for $J \subseteq I: H_J = \cap_{j \in J}H_j$. Assume that $\langle H_{I\setminus \{j\}} \mid j \in J \rangle = H_{I\setminus J}$ for all $J\subseteq I$. Then $\CC(G;(H_i)_{i \in I})$ corresponds to the chamber complex $\mathcal{C}(G, \cap_{i \in I} H_i, H_{I\setminus \{i\}}, i \in I)$ in a natural way:
    Given a coset complex $X =\CC(G;(H_i)_{i \in I})$ we consider as chambers $\mathcal{C}$ the maximal simplices of $X$ which are of the form $\{gH_{i} \mid i\in I\}$. Note that two maximal simplices $\{gH_i \mid i \in I \}, \{hH_i \mid i\in I\}$ are equal if and only if $h^{-1}g \in \cap_{i \in I} H_i$. Hence the set of maximal simplices corresponds to $G/H_I$. We say that $\{gH_i \mid i\in I\} \sim_j \{hH_i \mid i\in I\}$ if $gH_i = hH_i$ for all $i \in I\setminus \{j\}$ which is equivalent to $g H_{I\setminus \{j\}} = hH_{I\setminus\{j\}}$. 

    On the other hand, let $\mathcal{C}(G, \cap_{i \in I} H_i, H_{I\setminus \{i\}}, i \in I)$ be a chamber complex. For $J \subseteq I$ we denote for two chambers $C,D \in \mathcal{C}$ that $C \sim_J D$ if and only if there is a gallery from $C$ to $D$ with type set contained in $J$. For $C \in \mathcal{C}$ the $J$-residue of $C$ is the set $R_J(C) = \{D \in \mathcal{C} \mid C \sim_J D\}$. In our case, we have $$R_J(gH_I) = \{hH_I  \mid h \in g  \langle H_{I\setminus \{j\}} \mid j \in J \rangle \} = \{hH_I \mid h \in g  H_{I\setminus J}\}.$$ In particular, $R_{I\setminus \{i\}}(gH_I) = \{hH_I \mid h \in gH_j\}$. We define the set of vertices $X(0) := \{(R_{I\setminus \{j\}}(gH_I), j) \mid j \in J, g \in G\}$ which by the above argument is isomorphic to $\sqcup_{j \in J} G/H_j$. Furthermore, we say that two vertices $(R_{I\setminus \{j\}}(gH_I), j),(R_{I\setminus \{i\}}(hH_I), i)$ are adjacent if $i \neq j$ and $R_{I\setminus \{j\}}(gH_I) \cap R_{I\setminus \{i\}}(hH_I) \neq \varnothing$, where the latter is equivalent to $gH_j \cap hH_i \neq \varnothing$. From here, we build a simplical complex by defining that a set of vertices forms a simplex if and only if all vertices are pairwise adjacent.     
\end{remark}

\begin{definition}{\cite{TitsAmalgamees}} \label{def: direct limit}
    Let $(G_i)_{i \in I}$ be a family of groups and let $\phi_{ij}: G_i \to G_j$ for $(i , j) \in J\subseteq I \times I$ be a family of homomorphisms.

    The direct limit of the system $((G_i)_{i \in I}, (\phi_{ij})_{(i,j) \in J})$ is a group $G$ together with homomorphisms $f_i: G_i\to G$ such that $f_i = f_j \circ \phi_{ij}$ for all $(i,j) \in J$ and which satisfies the following universal property: 

    Given a group $H$ and homomorphisms $g_i: G_i \to H$ such that $g_i = g_j \circ \phi_{ij}$ for all $(i,j) \in J$ then there exists a unique homomorphism $\alpha: G \to H$ such that $g_i = \alpha \circ g_i$ for all $i \in I$.
        
    The group $G$ can be constructed explicitly by taking the free product of the groups $G_i, i\in I$ and making the identifications $\phi_{ij}(h) = h$ for all $h \in G_i, \forall (i,j) \in J$.
\end{definition}

\begin{remark} \label{rmk: complex of groups}
Let $(I,\leq)$ be a partially ordered set and let $J = \{(i,j)\in I \times I \mid i \leq j\}$. Let $(G_i)_{i \in I}$ be a family of groups with injective homomorphisms $\phi_{ij}: G_i \to G_j$ for all $i \leq j \in I$ such that $\phi_{ij} = \phi_{kj} \circ \phi_{ik}$ whenever $i \leq k\leq j$. Then $((G_i)_{i \in I},(\phi_{ij})_{(i,j)\in J})$ is called a simple complex of groups, see \cite[Definition 12.11]{bridson2011metric}.
In this case, the direct limit is called the fundamental group of the complex of groups. 
\end{remark}

The following result relates the notion of simple-connectedness of certain chamber systems and direct limits. It was first stated by J. Tits in \cite{TitsAmalgamees}, we present it in the formulation of \cite[Proposition 6.5.2]{ScharlauBuildings}.

\begin{proposition}{\cite[Proposition 6.5.2]{ScharlauBuildings}}\label{prop: m-simply-connected tool}
Let $G$ be a group, $B$ and $P_i, i \in I$, subgroups of $G$ such that $B \subset P_i$ for all $i \in I$. Assume that $G$ is generated by the $P_i$. Set $P_J:=\left\langle P_j: j \in J\right\rangle$ for all nonempty subsets $J \subseteq I$, and $P_{\varnothing}:=B$. The following are equivalent, for a natural number $m \leqslant|I|$:
\begin{enumerate}
    \item the chamber system $\mathcal{C}\left(G,B, P_i, i \in I\right)$ is $m$-simply-connected;
    \item $G$ is the direct limit of the $P_J,|J| \leqslant m$ with respect to the natural inclusions.
\end{enumerate}
\end{proposition}

\begin{remark}
    In the case when $m=2$ and $\lvert I \rvert = 3$ the result of \cite[Theorem 2.4]{abels} together with \ref{prop: m-simply-connected tool} show that the chamber system $\mathcal{C}(G,\cap_{i \in I} H_i, H_i, i\in I)$ is 2-simply connected if and only if $\CC(G; (\langle H_j \mid i \in I \setminus \{j\} \rangle)_{i \in I})$ is 1-connected, i.e. the zeroth and first homotopy group of the geometric realisation of the coset complex are trivial. 
    The requirement $\lvert I \rvert = 3$ is necessary here to assure that $\{\{i\}, \{i,j\} \mid i, j \in I \} = \{I\setminus \{i\}, I\setminus \{i,j\} \mid i,j \in I \}.$    

    For $m\neq 2$, the notions are, to the best of our knowledge, not related. Note that if a chamber system is $m$-simply connected then it is also $m+1$-simply connected. On the other hand, if a simplical complex is $n$-connected then it is also $n-1$-connected. 
\end{remark}

\section{Kac--Moody--Steinberg groups}\label{section: KMS group and root subgroups}

 Kac--Moody--Steinberg groups are defined as fundamental groups of certain complexes of finite $p$-groups, where $p$ is a fixed prime number. The finite $p$-groups involved in the construction are the positive unipotent subgroups of basic Levi subgroups of spherical type in a 2-spherical Kac--Moody group over a finite field. In \Cref{subsec: KMS groups}, we will define the Kac--Moody--Steinberg groups. We briefly discus the action of Kac--Moody groups on (twin) building and deduce properties of the KMS-groups from that action in \Cref{subsec: action on buildings}. In \Cref{subsection: quotients of KMS groups}, we show that nice infinite families of finite quotients of KMS groups exist. In particular, we will explain why, under certain assumptions, Kac--Moody--Steinberg groups are residually finite.

\subsection{Construction and properties} \label{subsec: KMS groups}

\begin{definition}\label{def: GCM}
$A =(A_{i,j})_{i,j \in I} \in \operatorname{Mat}_n(\mathbb{Z})$ is called a \emph{generalized Cartan matrix} (GCM) if
\begin{enumerate}
    \item $A_{i,i} = 2$ for all $i$;
    \item $A_{i,j} \leq 0$ for all $i \neq j$;
    \item $A_{i,j}=0 \iff A_{j,i}=0$.
\end{enumerate}

Every GCM $A=(A_{i,j})_{i,j \in I}$ gives rise to a Dynkin diagram in the following way. As vertex set, we take the index set $I$ and two vertices $i,j$ are connected by $|A_{i,j}|$ edges if $A_{i,j}A_{j,i} \leq 4$ and $|A_{i,j}| \geq |A_{j,i}|$, and these edges are equipped with an arrow pointing towards $i$ if $|A_{i,j}|>1$. If $A_{i,j}A_{j,i}>4$, the vertices $i$ and $j$ are connected by a bold-faced labelled edge with the ordered pair of integers $|A_{i,j}|,|A_{i,j}|$.

A GCM is \emph{irreducible} if there exists no non-trivial partition $I = I_1 \cup I_2$ such that $A_{i_1,i_2} = 0$ for all $i_1 \in I_1, i_2 \in I_2$.
\end{definition}

\begin{definition}
Let $A$ be a $(d+1)\times (d+1)$ generalized Cartan matrix with index set $I=\{0,\dots,d\}$, and let $J \subseteq I$.
\begin{enumerate}[label=(\roman*)]

    \item  The subset  $J$ is called \emph{spherical} if $A_J = (A_{i,j})_{i,j \in J}$ is of spherical type, meaning that the associated Coxeter group (see e.g. \cite[Proposition 4.22]{marquis2018introduction}) is finite (see e.g. \cite[Chapter 6.4.1]{bourbaki2008lie}).
Given $n  \geq 2$, $A$ is \emph{$n$-spherical} if every subset $J \subseteq I$ of size $n$ is spherical. 

\item We denote by $Q_A$ the set of spherical subsets of $I$ associated to the generalized Cartan matrix $A$.

\item A generalized Cartan matrix $A$ is \emph{purely $n$-spherical} if every spherical subset $J \subseteq I$ is contained in a spherical subset of size $n$.  (In particular, no set of size $n+1$ is spherical for a purely $n$-spherical generalized Cartan matrix.)

\item To a generalized Cartan matrix $A$ we can associate the following sets:
\begin{itemize}
    \item a set of simple roots $\Pi = \{\alpha_i \mid i \in I\}$,
    \item a set of real roots $\Phi \subseteq \bigoplus_{i \in I} \Z \alpha_i$,
    \item two sets, one of positive and one of negative real roots $\Phi^+= \bigoplus_{i \in I}\N \alpha_i \cap \Phi, \Phi^- = - \Phi^+$ .
\end{itemize}   
Note that $\Phi = \Phi^+ \sqcup \Phi^-$. More details can be found e.g. in \cite[Chapter 3.5]{marquis2018introduction}.
\end{enumerate}

\end{definition}

We assume for the rest of \Cref{section: KMS group and root subgroups} that $A$ is 2-spherical. Let $k$ be a finite field. To avoid degenerate cases, we assume $|k| \geq 4$.

Let $G=\mathcal{G}_A(k)$ be the (minimal) Kac--Moody group of type $A$ over a field $k$. This group is the quotient of the free product of the root subgroups
$$U_{\alpha} = \langle u_{\alpha}(t); \  t\in k \mid u_{\alpha}(t)u_{\alpha} (s) = u_{\alpha}(t+s) \text{ for all }s,t\in k \rangle \cong (k, +)$$ 
for all real roots $\alpha  \in \Phi$ and the torus 
$$T_k = \langle r^{\alpha_i} ; \  r \in k^\times, i \in I \mid r^{\alpha_i}s^{\alpha_i} = (rs)^{\alpha_i} \text{ for all } r,s \in k^\times\rangle \cong (k^\times)^{|I|}$$ with respect to certain relations, similar to those of Chevalley groups, see \cite[Definition 7.47]{marquis2018introduction}. In particular, if $A$ is a classical Cartan matrix, then $\GAk$ coincides with the Chevalley group of type $A$ over $k$. For more details about Chevalley groups, see also \Cref{section: Chevalley }.

To each spherical subset $J \subset I$, we associate the subgroup $U_J = \langle U_j \ | \ j \in J \rangle \leq \GAk$, where $U_j:=U_{\alpha_j} = \{u_{\alpha_j}(\lambda) \mid \lambda \in k \}$ is the root subgroup of $G$ associated to the simple root $\alpha_j$. We set $U_\varnothing = \{1\}$. We will refer to subgroups of the form $U_J$ as \emph{local groups}. 
The tuple of groups $(U_J)_{J \in Q_A}$, together with their natural inclusions, defines a simple complex of groups over the (finite) poset $Q_A$ (\Cref{rmk: complex of groups}). We denote by $\mathcal{U}_A(k)$ the fundamental group of that complex of groups, and we call it the \emph{Kac--Moody--Steinberg group of type $A$ over $k$} or KMS group for short. Thus, the group $\uak$ is the direct limit of the local groups with respect to all possible inclusions (see \Cref{def: direct limit}): 
$$\uak = *_{J\in Q_A} U_J / (U_J \hookrightarrow U_K, J \subseteq K \in Q_A). $$
We again denote by $U_J$ the canonical image of $U_J$ in $\uak$ for each $J \in Q_A$.

Since the base field $k$ is finite, each $U_J$ is a finite $p$-group, where $p$ is the characteristic of $k$ and the KMS group $\mathcal{U}_A(k)$ is finitely presented.

A homomorphism $f: \mathcal{U}_A(k) \to G$ is said to be \emph{injective on the local groups}, if $f \big|_{U_J}$ is injective for every $J$ in $Q_A$. Note that since $\mathcal{U}_A(k)$ is, by definition, the fundamental group of the complex of groups $(U_J)_{J \in Q_A}$, any family of maps $f_J: U_J \to G$ for $J \in Q_A$ such that $f_K \big|_{U_J} = f_J$ for $J \subseteq K$ factors through $\mathcal{U}_A(k)$.

\begin{remark} \label{rmk: presentation of KMS}
For two roots $\alpha, \beta \in \Phi$ we write $]\alpha, \beta[_\N = \left\{n_1 \alpha + n_2 \beta \in \Phi \mid n_1,n_2 \in \N^* \right\}$ and $[\alpha,\beta]_\N = \left]\alpha, \beta\right[_\N \cup \left\{\alpha,\beta\right\}$.

We have the following abstract presentation for the KMS groups:
$$\mathcal{U}_A(k)= \Bigl\langle u_\beta(t) \text{ for } t \in k, \beta \in [\alpha_i,\alpha_j]_{\mathbb{N}}, i,j \in I \ | \ \mathcal{R} \Bigr\rangle$$
where the set of relations $\mathcal{R}$ is defined as:

\indent for all $\{\alpha,\beta \} \subseteq [\alpha_i, \alpha_j]_{\mathbb{N}}, \ i,j \in I , \ t,u\in k$: 
$$[u_\alpha(t),u_\beta(u)] = \prod_{\gamma = k\alpha + l\beta \in ]\alpha,\beta[_{\mathbb{N}}} u_\gamma \left(C_{k,l}^{\alpha \beta} t^k u^l\right).$$
The constants $C_{k,l}^{\alpha,\beta}$ are called \emph{structure constants} and their precise values can be found for example in \cite{carter1989simple}.

Since $A$ is 2-spherical, the subgroups $U_\beta$, for $ \beta \in \left[ \alpha_i, \alpha_j\right] $, are contained in the group generated by the root groups $U_{\alpha_i}, U_{\alpha_j}$ (see \cite[Proposition 7]{AbTwinBuildings}). This can be seen using the commutator formula, similar to \cite[Lemma 3.13]{o2022high}.
\end{remark}

\begin{example} If we set $k=\mathbb{F}_p$ and $A_{ij}=-1$ for all $i \neq j$, then $\mathcal{U}_A(k)$ has the following presentation:

$$ \mathcal{U}_A(k) = \langle x_1 , \dots , x_d \ | x_i^p , \ [[x_i,x_j],x_j] \ (i\neq j) \rangle.$$

This example appears in \cite[Proposition 7.4]{ershov2010property} as an example of a Golod-Shafarevich group with property $(T)$ provided $d \geq 6$ and $p > (d-1)^2$.
\end{example}

\begin{remark} \label{rmk : Uak to U^+}
    Let $U^+ = U^+_A(k)=\langle u_\alpha(\lambda) \mid \alpha \in \Phi^+, \lambda \in k \rangle \leq \mathcal{G}_A(k)$ be the unipotent subgroup of the positive standard Borel subgroup of the Kac--Moody group. For example, in $\operatorname{SL}_n(k)$, the subgroup $U^+$ corresponds to the subgroup of upper triangular matrices with ones on the diagonal.

    There is a homomorphism from the KMS group to $U^+$:
\begin{align*}\psi: \mathcal{U}_A(k) \to U^+: \qquad u_j(\lambda) \mapsto u_{\alpha_j}(\lambda). \end{align*}
This homomorphism is injective on the local groups $U_J$ for $J \in Q_A$. 
\end{remark}

\subsection{Action on (twin) buildings and its consequences} \label{subsec: action on buildings}
A Kac--Moody group $\GAk$ together with its root subgroups $U_\alpha$ and torus $T_k$ forms an RGD-system $(\GAk, (U_\alpha)_{\alpha \in \Phi}, T_k)$. This was first observed by J. Tits, more details can be found in \cite[Chapter 8.8]{AB2008Buildings} and \cite[Theorem 7.69]{marquis2018introduction}. The RGD-system gives rise to the following groups
\begin{itemize}
    \item $U^+ := \langle U_\alpha \mid \alpha \in \Phi^+ \rangle, U^- := \langle U_\alpha \mid \alpha \in \Phi^- \rangle$,
    \item positive and negative Borel subgroup $B_+ := \langle T_k, U^+ \rangle$, $B_- := \langle T_k, U^- \rangle$,
    \item standard parabolic subgroups for $J\subseteq I$: $P^\pm_J:=\langle B_\pm, U_{\mp\alpha_j}; j\in J \rangle.$
\end{itemize}
Associated to the RGD-system is a twin building, which is a pair of two buildings $\Delta^+, \Delta^-$ together with additional structure which includes an opposition relation between the chambers (maximal simplices) of $\Delta^+$ and $\Delta^-$ (see \cite[Chapter 8]{AB2008Buildings}). As chamber systems, the buildings can be described as 
\begin{align*}
    \mathcal{C}^+ = \mathcal{C}(\GAk, B_+, P^+_{\{i\}}; i \in I) ,   \qquad \mathcal{C}^- = \mathcal{C}(\GAk, B_-, P^-_{\{i\}}; i \in I).
\end{align*}
They also can be described using the terminology of coset complexes (compare to \Cref{rmk: coco vs ch sys})
\begin{align*}
    \Delta^+ = \CC(\GAk; (P^+_{I\setminus \{i \}})_{i \in I}),     \qquad  \Delta^- = \CC(\GAk; (P^-_{I\setminus \{i \}})_{i \in I}).
\end{align*} 
In both descriptions, chambers correspond to cosets in $G/B_\pm$. Two chambers $gB_+, hB_-$ are opposite if and only if $g^{-1}h \in B_+B_-$ (see \cite[Equation (6.25)]{AB2008Buildings}). Two lower dimensional simplices $\tau_1,\tau_2$ are called opposite if they are contained in chambers $ \sigma_1, \sigma_2$ such that $\sigma_1$ is opposite $\sigma_2$.

We call the chamber $1B_+$ the positive fundamental chamber. We denote by $\mathcal{C}_{\text{op}}^- \subseteq \mathcal{C}^-$ the set of chambers opposite $1B_+$, and analogously by $\Delta_{\text{op}}^- \subseteq \Delta^-$ the complex opposite $1B_+$.

\begin{lemma}{\cite[Corollary 8.32]{AB2008Buildings}} \label{lemma: transitive action of U+}
The group $U^+$ acts sharply transitively on the set $\mathcal{C}_{\text{op}}^-$ of chambers opposite the positive fundamental chamber.
\end{lemma}
Let $\mathcal{U}_J := \langle U_\alpha \mid \alpha \in \Phi^+ \cap \oplus_{j \in J} \Z \alpha_j \rangle$. If the GCM $A$ is 2-spherical and $\lvert k \rvert \geq 4$ this group is equal to $U_J= \langle U_{\alpha_j} \mid j \in J \rangle$, which follows from the same argument as \Cref{prop: U_A(k) = U^+ part (a)}.
The following lemma is well known to experts, but we could not find a reference for the statement, hence we give a short proof.
\begin{lemma}\label{lemma: U intersect P}
    Let $\varnothing \neq J \subsetneq I$ be a proper subset of $I$. Then
    \begin{align*}
        U^+ \cap P_J^- = U_J.
    \end{align*}
\end{lemma}
\begin{proof}
    It is immediate from the definition that $\mathcal{U}_J \subseteq U^+ \cap P_J^-$.
    
To see the other inclusion, let $c_-= 1B_-$ be the negative fundamental chamber and consider its $J$-residue in $\mathcal{C}^-$.
\begin{align*}
    R_J(c_-) = \{ h B_- \mid h \in \langle P^-_{\{j\}} \mid j \in J \rangle \} = \{hB_- \mid h \in P^-_J\}.
\end{align*}
Its stabilizer in $G=\GAk$ is $\stab_G(R_J(c_-))=P_J^-$.

It follows from \cite[Corollary 5.30]{AB2008Buildings} that $R_J(c_-)=\mathcal{C}(P_J^-, B_-, P_j, j \in J)=:\mathcal{C}_J^-$ is a building in its own rights, of type $A_J:=(A_{ij})_{i,j \in J}$. Similarly, we get the building $\mathcal{C}_J^+:=\mathcal{C}(P_J^+, B^+, P_j^+, j\in J)$. This pair gives again rise to a twin building such that
$\mathcal{C}_{J,\text{op}} ^ - \cong R_J(c_-) \cap \mathcal{C}_{\text{op}}^-$.

Further note that the group $\mathcal{U}_J$ is the $U^+$ for the smaller twin building of type $A_J$.

Now, let $g \in U^+ \cap P_J^-$. Then $g$ stabilizes $R_J(c_-)$ by the above results and thus acts on $\mathcal{C}_J^-$. Furthermore, $gc_-$ is still opposite $c_+$ since $g \in U^+$. Hence $gc_- \in R_J(c_-) \cap \mathcal{C}_{\text{op}}^-$. By \Cref{lemma: transitive action of U+} and the above considerations, $\mathcal{U}_J$ acts transitively on $R_J(c_-) \cap \mathcal{C}_{\text{op}}^-$. Thus we can find $u \in \mathcal{U}_J=U_J$ such that $ugc_-=c_-$.  Hence $ug \in U^+ \cap \stab_G(c_-) =U^+ \cap B_- = \{1\}$, the last equality follows from \cite[Section 8.7]{AB2008Buildings}. Therefore we conclude $g \in U_J$ as desired.

\end{proof}

Using \Cref{rmk: converse of lemma} and \Cref{lemma: U intersect P}, we want to describe $\Delta_{\text{op}}^-$ as a coset complex. The chamber $1B_- = \{1P^-_{I\setminus \{i\}} \mid i \in I\}$ is clearly opposite $1B_+$. Thus it suffices to study the stabilizers inside $U^+$ of its vertices:
$$\stab_{U^+}(P^-_{I \setminus \{i\}}) = U^+ \cap P^-_{I\setminus \{i\}} = U_{I\setminus \{i\}}$$
where the last equality follows from \Cref{lemma: U intersect P}.
Hence $\Delta_{\text{op}}^- \cong \CC(U^+; (U_{I\setminus \{i\}})_{i\in I})$ as simplical complexes and $\mathcal{C}_{\text{op}}^- \cong \mathcal{C}(U^+, \{1\},U_{\alpha_i}; i \in I) $, where we used that $\cap_{i \in I} U_{\alpha_i}=\{1\}$.

In case the GCM that we started with is spherical, then the building $\Delta^+$ is finite and has an opposition relation within itself. In that case, two chambers $gB_+, hB_+$ are opposite if and only if $h^{-1}g \in B_+w_0,B_+$, where $w_0$ is the longest element in the Weyl group corresponding to the GCM. We analogously define $\Delta_{\text{op}} \subseteq \Delta^+$ to be the complex opposite the fundamental chamber $1B_+$. Again we have that 
$$\CC(U^+;(U_{I\setminus\{i\}})_{i \in I}) \cong \Delta_{\text{op}}.$$
This can be seen by a similar argument as above, but is also mentioned in \cite[Page 67]{AbTwinBuildings}.

In the rest of this section, we investigate how the above observations imply properties of the KMS-groups.

The following \Cref{prop: U_A(k)=U+} is well-known to experts, but the proof is scattered over different sources. For the convenience of the reader, we gather here all the pieces. Note that part (b) is a direct analogue to \cite[Corollary 1.2]{DevillersMühlherr2007}. In this formulation, it was first mentioned in \cite[Theorem 6.2]{ershov2008Golod}.

\begin{comment}
For the proof, we will need the following result, first stated by J. Tits in \cite{TitsAmalgamees}, in the formulation of \cite[Proposition 6.5.2]{ScharlauBuildings}. We call a chamber system $\mathcal{C}$ over $I$ $m$-simply-connected if it is connected (i.e. each two chambers can be connected by a gallery) and if each closed gallery is null-$m$-homotopic, see \cite[Chapter 6.1]{ScharlauBuildings} or \cite{DevillersMühlherr2007}.

\begin{proposition}(\cite[Proposition 6.5.2]{ScharlauBuildings})\label{prop: m-simply-connected tool}
Let $G$ be a group, $B$ and $P_i, i \in I$, subgroups of $G$ such that $B \subset P_i$ for all $i \in I$. Assume that $G$ is generated by the $P_i$. Set $P_J:=\left\langle P_j: j \in J\right\rangle$ for all nonempty subsets $J \subseteq I$, and $P_{\varnothing}:=B$. The following are equivalent, for a natural number $m \leqslant|I|$:
\begin{enumerate}
    \item the chamber system $\mathcal{C}\left(G,B, P_i, i \in I\right)$ is $m$-simply-connected;
    \item $G$ is the direct limit of the $P_J,|J| \leqslant m$ with respect to the natural inclusions.
\end{enumerate}
\end{proposition}
\end{comment}

\begin{proposition}\label{prop: U_A(k)=U+}
Let $k$ be an arbitrary field. Let $\psi: \mathcal{U}_A(k) \to U^+: u_j(\lambda) \mapsto u_{\alpha_j}(\lambda)$ be as in \Cref{rmk : Uak to U^+}.
\begin{enumerate}
    \item[(a)] If $A$ is 2-spherical and $|k| \geq 4$, then the map $\psi$ is surjective. \label{prop: U_A(k) = U^+ part (a)}
    \item[(b)]  \label{prop: U_A(k) = U^+ part (b)}   If $A$ is 3-spherical and $|k|\geq 5$, then $U^+$ is the direct limit of the root subgroups of rank at most 2 and thus
    $$\uak \cong U^+ $$
\end{enumerate}

\end{proposition}
\begin{proof}[Sketch of the proof]
\begin{comment}
To each Kac--Moody group, we can associate a twin BN-pair and thus a twin building, see for example \cite[Chapter 8]{AB2008Buildings} and \cite{marquis2018introduction} for more details. We consider the set $\Cop$ of opposite chambers to the positive fundamental chamber. Notice that the group $U^+$ acts sharply transitively on this set. We will consider the buildings as chamber systems. The set $\Cop$ can be described by the chamber system $\mathcal{C}(U^+, \{1\}, (U_{\alpha_i})_{i \in I})$ (see \cite[Chapter 1.4]{ScharlauBuildings} for more details about this notation) which has as chambers the elements of $U^+$ (or more precisely, the cosets $U^+/\{1\}$) and two chambers $c,c'\in U^+$ are $i$-adjacent, $i\in I$, if and only if $c^{-1}c' \in U_{\alpha_i}$. \\ 
\end{comment}

For part (a), note that $\Cop$ is gallery-connected if and only if $U^+ = \langle U_{\alpha_i} \mid i \in I \rangle$. The fact that $\Cop$ is gallery-connected under the assumption that $|k|\geq 4$ and that $A$ is 2-spherical can be found in \cite{ABRAMENKO1Muehl}. Since each $U_{\alpha_i}$ is in the image of $\psi$, this implies surjectivity.\\

For part (b), we want to apply \Cref{prop: m-simply-connected tool} for $m=2$ to $\Cop$. Thus, we need to show that $\Cop$ is simply 2-connected. \cite[Theorem 1.1]{DevillersMühlherr2007} reduces the problem to the following 
\begin{itemize}[label=-]
    \item for each rank 3 residue $R$ containing the positive fundamental chamber $c$, the set of chambers opposite of $c$ in $R$ is simply 2-connected;
    \item for each rank 2 residue $R$ containing the positive fundamental chamber $c$, the set of chambers opposite of $c$ in $R$ is connected.
\end{itemize}
Since we assume that $A$ is 3-spherical, each rank 3 and rank 2 residue will itself be contained in a spherical building. 

The required results for the rank 3 and rank 2 residues can be found in \cite{HornSimpleCon}, generalizing results of \cite{AbTwinBuildings}.\end{proof}

In the remainder of this section, we mention some properties of the subgroups $U_J$ of $\mathcal{U}_A(k)$. In particular, we describe their intersections as we will need them later on.

\begin{proposition} \label{prop: intersection of root subgroups}
    Assume that $A$ is 2-spherical and $k \geq 2$. Let $J, K \in Q_A$ with corresponding local groups $U_J, U_{K}$ then
    $$U_J \cap U_{K} = U_{J\cap K}.$$
\end{proposition}
\begin{proof}
By \cite[Theorem 4.3.8]{ScharlauBuildings}, we have for any $J, K \subseteq I$ that $P^-_J \cap P^-_K = P_{J\cap K}$. Combining this with \Cref{lemma: U intersect P}, we get
$$U_J \cap U_K = (U^+ \cap P_J^-) \cap (U^+ \cap P_K^-)= U^+ \cap (P_J^- \cap P_K^-) = U^+ \cap P_{J\cap K}^- = U_{J\cap K}.$$
\end{proof}

\subsection{Finite quotients of Kac--Moody--Steinberg groups}\label{subsection: quotients of KMS groups}
In this section, we describe how to construct infinite families of finite quotients of Kac--Moody--Steinberg groups. 
To do so, we rely on the properties of the natural map $\psi : \uak \to U^+$ and the fact that $U^+$ is residually-$p$. The local injectivity of the maps we consider allows us to conclude.
The details of those arguments can be found in this section and the explicit construction of infinite families of finite quotients of the KMS group as well.

We will use the following observation.
\begin{observation} \label{observation: equivalence to stupid property}
    Let $f:X \to Y$ be a homomorphism between two groups $X,Y$ and let $A,B \leq X$ be subgroups such that $f|_A$ and $f|_B$ are injective. Then 
    \begin{align*}
        f(A)\cap f(B) = f(A\cap B) \iff \ker(f) \cap \bigl[ (A\setminus (A\cap B)) \cdot (B \setminus (A \cap B)) \bigr ]= \varnothing
    \end{align*}
    where for any sets $A',B' \subseteq X$ we define $A' \cdot B' := \{a \cdot b \mid a \in A', b \in B'\}$.
\end{observation}

\begin{lemma}\label{lemma: composition and stupid property}
Let $\psi: \uak \to U^+$ be as in \Cref{rmk : Uak to U^+} and write $\mathcal{U}_J = \psi(U_J)$ for all $J \in Q_A$. Let $H$ be any group and $\phi:U^+ \to H$ a homomorphism satisfying 
    \begin{enumerate}
        \item $\phi|_{\mathcal{U}_J}$ is injective for all $J \in Q_A$;
        \item $\phi(\mathcal{U}_J) \cap \phi(\mathcal{U}_K) = \phi(\mathcal{U}_J \cap \mathcal{U}_K)$ for all $J,K \in Q_A$.
    \end{enumerate}
    Then the composition $f = \phi \circ \psi: \uak \to H$ satisfies, for all $J,K \in Q_A$
    \begin{align*}
        f(U_J) \cap f(U_K) = f(U_J \cap U_K)
    \end{align*}
    and $f$ is injective on the local groups.
\end{lemma}
\begin{proof}
The map $f$ is injective on the local groups as it is the composition of two maps that are injective on the local groups.

Let $J,K \in Q_A$ and set $U_J'= U_J \setminus (U_J \cap U_K)$ and $ \mathcal{U}_J'= \mathcal{U}_J \setminus (\mathcal{U}_J \cap \mathcal{U}_K)$ and similarly for $U_K'$ and $\mathcal{U}_K'$. Using \Cref{observation: equivalence to stupid property} it suffices to show $\ker(f) \cap U_J'\cdot U_K' = \varnothing$. Note that $\ker(f) = \ker(\phi \circ \psi) = \psi^{-1}( \ker(\phi))$. Moreover,
\begin{equation}\label{equivalence : intersection root groups}
    \psi^{-1}(\ker(\phi)) \cap U'_J \cdot U_K' = \varnothing \iff \ker(\phi ) \cap \psi(U'_J \cdot U'_K) = \varnothing.
\end{equation}
Since $\psi$ is a homomorphism injective on the local groups, we have $\psi(U'_J \cdot U_K') = \mathcal{U}_J' \cdot \mathcal{U}_K'$. By assumption, the equality on  the right-hand side of the equivalence in Equation \ref{equivalence : intersection root groups} holds.
\end{proof}

\begin{theorem} \label{thm: residually finite}
    Let $A$ be a 2-spherical GCM and let $k$ be a finite field of characteristic $p$. Then the group $U^+$ is residually-$p$. If $A$ is 3-spherical and $|k| \geq 5$, this also holds for $\uak$.
\end{theorem}
\begin{proof}
    This is a consequence of \cite[1.C Lemma 1]{remy2006topological}. For the readers convenience, we provide an argument for the residually finiteness of $U^+$.
    
    We have that $U^+$ is contained in the stabilizer of the fundamental chamber in the positive building associated to the Kac--Moody group. Thus, the action of $U^+$ stabilizes balls of radius $r \in \N$ around the fundamental chamber (in the chamber complex). The point-wise fixers of these balls are then finite index normal subgroups of $U^+$. Since $U^+$ acts faithfully on the building, the conclusion follows.
    In the case where $A$ is 3-spherical and $|k| \geq 5 $, the groups $U^+$ and $\uak$ are isomorphic by \Cref{prop: U_A(k)=U+}.
\end{proof}

Note that the sets $\mathcal{U}_J \cdot \mathcal{U}_K$ are finite since they can be embedded into the product of two finite sets.\\
\Cref{thm: residually finite} above together with the fact that $|\bigcup_{J,K \in Q_A} \mathcal{U}_J\cdot \mathcal{U}_K| < \infty$  and \Cref{lemma: composition and stupid property} implies the following corollary.

\begin{corollary} \label{cor: inf fam of f.i. subgrps}
Let $A$ be a 2-spherical GCM and let $k$ be a finite field of characteristic $p$ such that $|k| \geq 4$. Then, there is an infinite family of finite index normal subgroups $N_1,N_2,\dots  \unlhd_{\text{f.i.}} U^+$ such that
\begin{enumerate}
    \item $N_i \cap \left( \bigcup_{J,K \in Q_A} \mathcal{U}_J\cdot \mathcal{U}_K \right) = \{1\}$ \label{cor: inf fam of f.i. subgrps condition 1} for all $i$, and;
    \item $\lim_{i\to \infty} |U^+/N_i| = \infty$.
\end{enumerate}
Moreover, the quotient maps 
$$ \varphi_i = \phi_i \circ \psi:\uak \overset{\psi}{\longrightarrow} U^{+} \overset{\phi_i}{\longrightarrow} U^{+}/N_i $$ are injective on the local groups and for any $J, K \in Q_A$ we have
$\varphi_i(U_J) \cap \varphi_i(U_{K}) = \varphi_i(U_J \cap U_{K})$.

\end{corollary}

\section{Main theorem}\label{section: main theorem}
In this section, we give the proof of the main theorem. We will first show results about the connectivity of the links, and then about the degree of the coset complex we construct. After doing so we can prove the expansion properties of the links and finally describe some symmetries of our construction.

For this section, we fix some notation. Let $(A_{i,j})_{i,j\in I}, I = \{0,\dots,d\}$ be a non-spherical generalized Cartan matrix of rank $d+1$ that is $d$-spherical and purely $d$-spherical, i.e., the set of spherical subsets is $Q_A = \mathcal{P}(I) \setminus \{I\}$. Furthermore, let $k$ be a finite field of order $|k| = q=p^m \geq 4$ and let $U_i \subseteq \mathcal{G}_A(k)$ be the root subgroup associated to the simple root indexed by $i \in I$ in the Kac--Moody group $\mathcal{G}_A(k)$. Let $G$ be a finite group, $\phi:  \uak \to G$ be a homomorphism from the KMS group to $G$ that is injective on the local groups and such that $\langle \phi(U_0), \dots, \phi(U_d) \rangle = G$ (i.e. $\phi$ is surjective). We further require that $\phi$ satisfies the following intersection property:

\begin{equation}\label{essential property}
\tag{IP}
\phi(U_J) \cap \phi(U_{K}) = \phi(U_J \cap U_{K}) \text{ for all } J, K \in Q_A.
\end{equation}

We introduce the following notation.
\begin{definition}
For $i \in I$, we define $\hi = I \setminus \{i\}$ and
$$H_i:= \phi(U_{\hi}).$$
For $\varnothing \neq T \subseteq I$, let 
$$H_T:= \bigcap_{i \in T} H_i.$$
We set $H_{\varnothing} = \langle H_0,\dots,H_d \rangle$.
\end{definition}

The following lemma is useful when working with the images of the local groups and allows us to describe such images in terms of subsets of the index set $I$.

\begin{lemma}\label{lemma: relation H and U}
For $\varnothing \neq T \subseteq I$ we have $H_T = \phi(U_{I\setminus T})$.
\end{lemma}
\begin{proof}
We have:
\begin{align*}
    H_T = \bigcap_{i \in T} H_i = \bigcap_{i \in T} \phi\left(U_{\hi} \right) \overset{(*)}{=} \phi\left(\bigcap_{i \in T} U_{\hi}\right) \overset{(**)}{=} \phi\left(U_{\bigcap_{i \in T} I\setminus \{i\}}\right) = \phi\left(U_{I\setminus T}\right) 
\end{align*}
where $(*)$ follows from \Cref{essential property} and $(**)$ from \Cref{prop: intersection of root subgroups}.

\end{proof}

In the following theorem, the main result of this paper, we prove that a well-chosen coset complex associated with a generalized Cartan matrix $A$, a field $k$ and a complex of groups $(U_J)_{J\in Q_A}$ is a $\lambda$-spectral HDX. We maintain the notation used above in the following theorem.

\begin{theorem}\label{thm: main thm}
Let $G, H_i, i \in I$ be as above and consider the coset complex
    $$X = \CC(G; (H_i)_{i\in I}).$$
Then, $X$ is a pure $d$-dimensional simplicial complex and 
\begin{enumerate}
    \item \label{main thm: item connected} for any $\sigma \in X(i), i \leq d-2$ we have that $\lk_X(\sigma)$ is connected. In particular, $X$ is connected.
    \item \label{main thm: item bdd degree} There exists $\Delta \in \N$ such that for all $v \in X(0)$: $|\{\sigma \in X(d): v \subseteq \sigma \}| \leq \Delta$. 
    \item \label{main thm: item local expansion} For any $\sigma \in X(d-2): \lambda_2(\lk_X(\sigma)) \leq \gamma$, where $\gamma$ is specified below.
    \item \label{main thm: item nice action} $G$ acts sharply transitively on the maximal faces of $X$.
\end{enumerate}

In particular, if $k$ and $d$ are chosen such that $\gamma \leq \frac{1}{d}$, $X$ is a $\gamma'$-spectral HDX for $\gamma'= \frac{\gamma}{1-(d-1)\gamma}$. The exact values for $\gamma$ are: 
\begin{itemize}
    \item if the Dynkin diagram associated to $A$ has at most single edges, we have $\gamma = \frac{1}{\sqrt{q}}$;
    \item if the Dynkin diagram associated to $A$ has at most double edges, we have $\gamma = \sqrt{\frac{2}{q}}$;
    \item if the Dynkin diagram associated to $A$ has a triple edge, we have $\gamma = \sqrt{\frac{3}{q}}$. 
\end{itemize}

\end{theorem}

The following corollary describes how to construct a family of HDX from a Kac--Moody--Steinberg group.

\begin{corollary}\label{corol: main corollary}
    Let $d \geq 3$, $A$ a 2-spherical GCM, and $k$ a finite field of size at least 4. 
    By \Cref{cor: inf fam of f.i. subgrps}, we have an infinite family of quotient maps $\varphi_i: \uak \to G_i$, where $G_i = U^+/N_i$ as in the notation of the corollary. The maps $\varphi_i$ are injective on the local groups and satisfy \Cref{essential property}. Let $H_{j}^{(i)}:= \varphi_i(U_{\hj}), j \in I$ be the images of the local groups and set $X_i:= \CC(G_i, (H_j^{(i)})_{j\in I})$.

    Then $\left( X_i \right)_{i \in \N}$ is a family of bounded degree high-dimensional spectral expanders.
\end{corollary}

\begin{remark}
The \emph{basic construction} given in \cite[Chapter II.12.11]{bridson2011metric} using $(U_J)_{J \in Q_A}$, together with their natural inclusions, as the complex of groups yields a simplicial complex of dimension $d$, where $d$ is the size of the maximal spherical subsets of $I$. Our goal was initially to explore the expansion properties of this simplicial complex, obtained via maps $\phi: \mathcal{U}_A(k) \to G$, where $G$ is a finite group and such that the map $\phi$ is injective on the local groups $U_J$ for each $J \in Q_A$.

However, if $\phi: \mathcal{U}_A(k) \to G$ satisfies \Cref{essential property}, the simplicial complex obtained using the construction described in  \cite[Chapter II.12.11]{bridson2011metric} can be simply described in terms of coset complexes and thus is amenable to an approach similar to \cite{o2022high}, as is done in this paper.

It is worth noting that the basic construction in \cite[Chapter II.12.11]{bridson2011metric} yields another natural way to associate simplicial complexes to complexes of groups and should therefore be kept in mind as a potential source for other high-dimensional expanders.
\end{remark}

\vspace{0.5cm}
The following sections are devoted to the proof of \Cref{thm: main thm}. In \Cref{subsection: connectivity links} we prove that all the links of $X$ are connected, then in \Cref{subsection: bounded degree } we bound the number of maximal faces containing a given vertex of $X$. In \Cref{subsection: link expansion} we prove spectral expansion for the $(d-2)$-dimensional links of $X$ and finally in \Cref{subsection: symmetries} we describe some symmetries of our complex. All these sections combined prove the different statements of \Cref{thm: main thm}.

\subsection{Connectivity of the links}\label{subsection: connectivity links}

\begin{proposition} \label{prop: strucutre of links}
    Let $\sigma \in X(i)$ for $ i \leq d-2$ be of type $\varnothing \neq T \subseteq I$. Then 
    $$\lk_X(\sigma) = \CC\left(H_T ,\left(H_{T \cup\{i\}}: i \notin T\right)\right) \cong \CC\left(U_{I \setminus T}, \left(U_{I \setminus (T \cup \{i\})}\right)_{i \in I \setminus T}\right). $$
\end{proposition}
\begin{proof}
    This follows directly from the fact that $\phi |_{U_{I \setminus T}}: U_{I\setminus T} \to \phi(U_{I \setminus T}) = H_T$ is an isomorphism of groups that maps $H_{T \cup \{i\}}$ to $U_{I \setminus (T \cup \{i\})}$ (see \Cref{lemma: relation H and U}).
\end{proof}

\begin{proposition} \label{prop: generating root subgroups from others}
For $\varnothing \neq T \subsetneq I$ we have
    $$U_{I\setminus T} = \langle U_{I \setminus (T \cup \{i\})}: i \in I\setminus T \rangle .$$
Furthermore, we have  
$$H_\varnothing = G.$$
\end{proposition}
\begin{proof}
Let $\varnothing \neq T \subsetneq I$. Then
\begin{align*}
     \langle U_{I \setminus (T \cup \{i\})}: i \in I\setminus T \rangle &=  \langle U_j: j \in I \setminus (T \cup \{i\}) , i \in I \setminus T \rangle \\
     &= \langle U_j: j \in \cup_{i \in I \setminus T} I \setminus (T \cup \{i\}) \rangle \\
     &= \langle U_j: j \in I \setminus T \rangle = U_{I\setminus T}.
\end{align*}

For $T = \varnothing$ we get:
$$G \supseteq H_\varnothing = \langle H_0,\dots,H_d \rangle \supseteq \langle \phi(U_0),\dots, \phi(U_d) \rangle = G$$
where the last equality follows from the surjectivity assumption on $\phi$.
\end{proof}

   From \Cref{prop: connectivity of CC} together with \Cref{prop: strucutre of links} and \Cref{prop: generating root subgroups from others} we get that for every $\sigma \in X(i), i \leq d-2$ the link $\lk_X(\sigma)$ is connected.

\subsection{Bounded degree}\label{subsection: bounded degree }

We want to consider the number of maximal faces that contain a given arbitrary vertex $gH_i$.

Maximal faces in $X$ are of the form $\{hH_0,\dots,hH_d\}$. Two maximal faces $\{hH_0,\dots,hH_d\}$, $\{h'H_0,\dots,h'H_d\}$ coincide if and only if $h^{-1}h' \in \bigcap_{i \in I} H_i = H_I = \phi(U_\varnothing) = \{1\}$. Thus, every element in $G$ gives rise to a unique maximal face.

Which maximal faces contain a given vertex $gH_i$? Note that:
$$gH_i \in \{hH_0,\dots,hH_d\} \iff gH_i = hH_i \iff h \in gH_i.$$
Thus
$$\deg(gH_i) = |H_i|=|U_{\hi}| = |\langle U_j \mid j \in I \setminus \{i\}\rangle |.$$

The cardinality of $ U_{\hi}$ is equal to $|k|^\ell$ where $\ell$ is the number of positive roots in the root system of type $(A_{m,j})_{m,j \in I\setminus \{i\}}$ by \cite[Theorem 5.3.3]{carter1989simple}. The number of positive roots in the root system of type $(A_{m,j})_{m,j \in I\setminus \{i\}}$ is bounded from above by $(|I|-1)^2 = d^2$ if $(A_{m,j})_{m,j \in I\setminus \{i\}}$ is of classical type (i.e. $A_d,B_d,C_d,D_d$). Exact numbers can be found in \cite[Chapter 6, Section 4]{bourbaki2008lie}.

Thus the coset complexes we obtain are of bounded degree and the bound depends only on the choice of the generalized Cartan matrix and the field we chose for the Kac--Moody--Steinberg group $\uak$, proving \Cref{main thm: item bdd degree} of \Cref{thm: main thm}.

\subsection{Link expansion}\label{subsection: link expansion}
Let $\sigma \in X(d-2)$ of type $T\subseteq I$. Then $|T| = d-1$ and thus $I \setminus T = \{i,j\}$ has two elements. From 
\Cref{prop: strucutre of links} we get that
$$\lk_X(\sigma) \cong \CC(U_{i,j}; U_i,U_j).$$

The idea for proving the expansion of the graphs $\CC(U_{i,j};U_i,U_j)$ is by considering them as $\lvert k \rvert$-regular subgraphs of the spherical building of the corresponding type. This is indeed possible since $U_{i,j}$ is the positive unipotent subgroup in the Chevalley group of type $A_{\{i,j\}}$ and using \Cref{subsec: action on buildings}. The corresponding spherical building is $\lvert k \rvert + 1$-regular and its spectrum is known. Hence the following lemma gives the desired bound on the second smallest eigenvalue of the random walk matrix of $\CC(U_{i,j};U_i,U_j)$. The proof of the lemma only uses basic linear algebra tools and can be found in \cite[Lemma 5.5]{harsha2019note}.

\begin{lemma}\label{lemma: spectrum of regular subgraph}
    Let $X$ be a $d$-regular, simple graph that is an induced subgraph of a $D$-regular graph $Y$. Then
    $$\lambda_2(X)\leq \frac{D}{d} \lambda_2(Y).$$
\end{lemma}

Thus, we investigate the spectrum of the spherical buildings corresponding to the Chevalley groups of type $A_2,B_2,G_2$ corresponding to the cases where $i,j$ are connected by a single, double or triple edge in the Dynkin diagram, respectively.

%------------- New version --------------------
The eigenvalues of the Laplacian matrix $\Delta = \id - M$ (where $M$ is the random walk matrix) of the spherical building of rank 2 are computed in \cite[Proposition 7.10]{garland1973p}. Note that in Garlands work we have $q' = q(1) = q(2) = \lvert k \rvert$ (since per definition $1+q(i) = \lvert P_{\{i \}} / B \rvert = \deg(1P_{\{i\}}) = \lvert k \rvert +1$). 
Hence we can deduce from \cite[Proposition 7.10]{garland1973p} the following values for the spectrum of the random walk matrix depending on the type.

\subsubsection*{The $A_2$ case}
Let $k = \mathbb{F}_q$ and let $M$ be the random walk matrix of the spherical building of type $A_2$ over $k$. Then the spectrum is
$$\operatorname{spec}(M) = \left\{1, \frac{\sqrt{q}}{q+1}, -\frac{\sqrt{q}}{q+1}, -1 \right\}$$
and hence 
$$\lambda_2(\CC(U_{i,j};U_i,U_j)) \leq \frac{q+1}{q} \cdot \frac{\sqrt{q}}{q+1}= \frac{1}{\sqrt{q}}.$$

\subsubsection*{The $B_2$ case}
Let $k = \mathbb{F}_q$ and let $M$ be the random walk matrix of the spherical building of type $B_2$ over $k$. Then the spectrum is $$\operatorname{spec}(M) = \{1, \frac{\sqrt{ 2q }}{q+1}, 0,- \frac{\sqrt{ 2q }}{q+1},-1 \}$$
and hence $$\lambda_{2}(\mathcal{CC}(U_{\left\{ i,j \right\}};U_{i},U_{j})) \leq \frac{(q+1)\sqrt{ 2q }}{q(q+1)}=\sqrt{ \frac{2}{q} }.$$

\subsubsection*{The $G_2$ case}
Let $k = \mathbb{F}_q$ and let $M$ be the random walk matrix of the spherical building of type $G_2$ over $k$. Then the spectrum is
$$\operatorname{spec}(M) = \{1, \frac{\sqrt{ 3q }}{q+1}, \frac{\sqrt{ q }}{q+1}, 0, -\frac{\sqrt{ q }}{q+1}, -\frac{\sqrt{ 3q }}{q+1},-1\}$$ 
and hence 
$$\lambda_{2}(\mathcal{CC}(U_{\left\{ i,j \right\}};U_{i},U_{j})) \leq \frac{(q+1)\sqrt{ 3q }}{q(q+1)}=\sqrt{ \frac{3}{q} }.$$

\begin{remark}
    If $q$ is prime, the bounds on $\lambda_2(\CC(U_{i,j};U_i,U_j)$ are sharp in the $A_2$ and $B_2$ case. This can be seen using the concept of representation angle and the results of \cite[Proposition 7.3]{PEC2022Hyperbolic}.
\end{remark}

\subsection{Symmetries of the construction}\label{subsection: symmetries}
Similar to the constructions of \cite{kaufman2018construction} and \cite{o2022high} but different to the first construction of HDX giving rise to Ramanujan complexes \cite{LubSamVishRamanujan}, our simplicial complexes are highly symmetric in the following sense.
\begin{proposition}
$G$ acts sharply transitively on the maximal faces of $X$.
\end{proposition}
\begin{proof}
Following \Cref{lem: action of G on CC}, it remains to check that $H_0 \cap \dots \cap H_d = \{1\}$. But this follows directly since $H_0 \cap \dots \cap H_d= \phi(U_{\widehat{0}} \cap \dots \cap U_{\widehat{d}}) = \phi(U_\varnothing)= \{1\}$.
\end{proof}

This proves \Cref{main thm: item nice action} of \Cref{thm: main thm} and thus finishes the proof.

\section{Quotients in Chevalley groups}\label{section: Chevalley }
In this chapter, we describe a connection between spherical Chevalley groups and affine Kac--Moody--Steinberg groups. This allows us to describe our construction with a very concrete example in $\operatorname{SL}_3$ (this is done in \Cref{subsection: quot in Chev: sl3}) and to compare it to previous constructions of high-dimensional expanders using coset complexes from \cite{kaufman2018construction} and \cite{o2022high} in \Cref{subsection: quot in chev: comparison}.

\subsection{Definitions and properties} \label{subsection: quot in Chev: Def and Prelim}
We start by recalling the definition of Chevalley groups, the connection between spherical and affine root systems and how this leads to a map from certain KMS groups to Chevalley groups over polynomial rings.

The most general way to define Chevalley groups is by introducing the \emph{Chevalley--Demazure group scheme}, which is a functor from the category of commutative unital rings to the category of groups. This is done by considering certain Hopf-algebras over $\Z$ that are related to semisimple Lie groups. The notion was established in the work of Chevalley \cite{ChevalleySchemasdegrpsSemiSimple} and Demazure \cite{DemazureSchemasengrpsReductifs}, a formal definition can be found in \cite{Abe1969ChevalleySchemes}.

Over fields, Chevalley groups can be defined by giving their so-called Steinberg presentation, see e.g. \cite[Theorem 12.1.1]{carter1989simple}. Note that some definitions consider the group defined below modulo its center as Chevalley groups. We will not take this into account here.

\begin{definition}
Corresponding to any irreducible, spherical Cartan matrix $\Ao$ with root system $\mathring \Phi$ of rank at least 2, and any finite field $k$, there is an associated universal (or simply connected) Chevalley group, denoted $\chev_{\Ao}(k)$. Abstractly, it is generated by symbols $x_\alpha(s)$ for $\alpha \in \mathring \Phi$ and $s \in k$, subject to the relations
$$
\begin{aligned}
x_\alpha(s) x_\alpha(u) & =x_\alpha(s+u) \\
{\left[x_\alpha(s), x_\beta(u)\right] } & =\prod_{i, j>0} x_{i \alpha+j \beta}\left(C_{i j}^{\alpha, \beta} s^i u^j\right) \quad(\text {for } \alpha+\beta \neq 0) \\
h_\alpha(s) h_\alpha(u) & =h_\alpha(s u) \quad(\text { for } s,u \neq 0), \\
\text { where } \quad h_\alpha(s) & =n_\alpha(s) n_\alpha(-1) \\
\text { and } \quad n_\alpha(s) & =x_\alpha(s) x_{-\alpha}\left(-s^{-1}\right) x_\alpha(s) .
\end{aligned}
$$
Note that the structure constants $C_{ij}^{\alpha, \beta}$ are the same as for KMS groups (see \Cref{rmk: presentation of KMS}).
\end{definition}
\begin{remark}
Let $k$ be a field and let $k[t]$ denote the polynomial ring in one variable over $k$. The simply-connected Chevalley group $\chev_{\Ao}(k[t])$ is, similar to the case of a Chevalley group over a field, generated by elements $x_\alpha(s)$ for $\alpha \in \mathring \Phi, s \in k[t]$ that satisfy the relations above, where for the third relation we have to add the extra assumption that $u,s$ are invertible in $k[t]$. This can be found e.g. in \cite[Chapters 7,9 and 13]{vavilov1996chevalley}.
\end{remark}

Let $\Ao = (A_{i,j})_{i,j \in \mathring{I}}, \mathring{I} = \{1,\dots d \}$ be an irreducible (spherical) Cartan matrix with spherical root system $\mathring \Phi$ and set of simple roots $\mathring \Pi = \{\alpha_1,\dots,\alpha_d\}$. This root system has a unique highest root $\gamma \in \mathring \Phi$, which means that $\gamma$ is such that for any $\alpha_i \in \mathring \Pi$ we have $\gamma + \alpha_i \notin \mathring \Phi$. 

Whenever we remove one root from the set $S:= \{-\gamma, \alpha_1,\dots,\alpha_d\}$ we get a set of simple roots (although they might not generate the full root system $\mathring \Phi$ but a subsystem).

We set $\alpha_0 = -\gamma$ and $I = \{0,\dots, d\}$. The idea now is to consider the set $\Pi = \{\abar_0,\dots,\abar_d\}$ as set of simple roots of a root system $\Phi$ such that for all $i,j \in I$ the roots $\abar_i,\abar_j$ generate the same rank two subsystem as $\alpha_i, \alpha_j$ in $\mathring{\Phi}$. The resulting root systems $\Phi$ with simple roots $\Pi = \{\abar_0,\dots,\abar_d \}$ and generalized Cartan matrix $A$ are described by the Dynkin diagrams in \cite[§4.8]{KacInfDimLieAlg} (or \cite[Table 5.1 Aff 1]{marquis2018introduction}) and are of affine (untwisted) type. Note that $A$ is of rank $d+1$, is $d$-spherical and purely $d$-spherical.

For a root $\beta = \sum_{i \in \mathring{I}} \lambda_i \alpha_i \in \mathring{\Phi}$, $\lambda_i \in \Z$, we write $\overline{\beta} = \sum_{i \in \mathring{I}} \lambda_i \abar_i \in \Phi$. We set $\overline{\delta} = \abar_0 + \overline{\gamma} \in \bigoplus_{i \in I} \Z \abar_i$ which is not in $\Phi$ but it is what is called an imaginary root of the root system associated to $A$. For our purposes, it suffices to think of it as an element in $\bigoplus_{i \in I} \Z \abar_i$. From \cite[Section 7.4]{KacInfDimLieAlg} it follows that $\Phi = \{ a_{\alpha,m}:= \abar + m \overline{\delta} \mid m \in \Z, \alpha \in \mathring \Phi \}$. Note that $\abar_0 = a_{-\gamma,1}$. A root $a_{\alpha,m}$ is in $\Phi^+$ if and only if $m\geq 1$ or $m=0$ and $\alpha \in \mathring \Phi^+$ (see e.g. \cite[proof of Theorem 7.90]{marquis2018introduction}).

Before we continue with the construction, we fix a few notations: $k[t]$ denotes the polynomial ring in one variable $t$ over $k$. For a polynomial $f \in k[t]$, let $(f) = \{g \cdot f \mid g \in k[t]\}$ denote the ideal generated by $f$ in $k[t]$ (e.g. $(t)$ is the set of polynomials without constant term) and let $kt^n = \{\lambda t^n \mid \lambda \in k\}$ be the set of all scalar multiples of $t^n$.

Let $k$ be a finite field of order $|k| \geq 4$ and assume that $\operatorname{rank}(\mathring \Phi) \geq 2$ (thus we have $\uak \twoheadrightarrow U^+$). \cite[Theorem 7.90]{marquis2018introduction} implies that the following is a well-defined injective homomorphism
\begin{alignat*}{1}
    \Tilde{\phi}: U^+ & \to \chev_{\Ao}(k[t]) = \langle x_\alpha(g) \mid \alpha \in \mathring \Phi, g \in k[t] \rangle \\
    u_{a_{\alpha,m}}(\lambda) & \mapsto x_\alpha(\lambda t^m).
\end{alignat*}

Let $f \in k[t]$ be an irreducible polynomial of degree $\ell \geq 2$. Then $k[t] = (t) + (f)$ and we denote by $ \pi_f: k[t] \to k[t]/(f)$ the projection to the quotient of the polynomial ring by the ideal generated by $f$. Thus, we have $\pi_f(k[t]) = \pi_f((t))$. The restriction of $\pi_f$ to $\{g \in k[t] \mid \deg(g) < \ell \}$ is injective, since $\{g \in k[t] \mid \deg(g) < \ell \} \cap (f) = \{0\}$. 

Using the functoriality of the Chevalley group functor, we can go from the Chevalley group over $k[t]$ to the one over $k[t]/(f)$ by applying $\pi_f$ "entry-wise". We furthermore precompose with the map $\psi: \uak \to U^+$ from \Cref{rmk : Uak to U^+}. We get the following well-defined map
\begin{alignat*}{3}
    \varphi_f: \uak & \overset{\psi}{\longrightarrow} U^+ && \overset{\Tilde{\phi}}{\longrightarrow} \chev_{\Ao}(k[t]) && \overset{\pi_f}{\longrightarrow} \chev_{\Ao}(k[t]/(f))\\
    u_i(\lambda) & \mapsto u_{\abar_i}(\lambda) && \mapsto \begin{cases} x_{\alpha_i}(\lambda)  & i \neq 0 \\ x_{-\gamma}(\lambda t) & i = 0 \end{cases} && \mapsto \begin{cases} x_{\alpha_i}(\lambda + (f))  & i \neq 0 \\ x_{-\gamma}(\lambda t + (f)) & i = 0 \end{cases}.
\end{alignat*}
We will use the notation $\phi = \Tilde{\phi} \circ \psi$. 

The goal is now to show that the map $\varphi_f$ satisfies the requirements of \Cref{thm: main thm}, i.e. local injectivity, surjectivity and that the image of the intersection of two local groups is the intersection of the images. This will be done in the next section. In \Cref{subsection: quot in Chev: sl3}, we describe our construction in the most basic case, namely the one where $\chev_{\Ao} = \operatorname{SL}_3$. In this case, we can describe explicitly the subgroups $H_i$ yielding a family of HDX.

\subsection{The general case} \label{subsection: quot in Chev: the general case}
We keep the notation from \Cref{subsection: quot in Chev: Def and Prelim} above.

\begin{definition}
 For $J \in Q_A= \{J \subseteq I \mid J \text{ spherical} \}$, note that $\{\Bar{\alpha}_j \mid j\in J \} \subset \Pi$ and $\{\alpha_j \mid j \in J\} \subset \mathring \Pi \cup \{\alpha_0\}$ are sets of linear independent roots that generate roots subsystems $\Phi_J = \bigoplus_{j \in J} \Z \abar_j \cap \Phi$ and $ \mathring \Phi_J = \bigoplus_{j \in J} \Z \alpha_j \cap \mathring \Phi$, respectively. Their set of positive roots is denoted by $\Phi_J^+ = \bigoplus_{j \in J} \N \abar_j \cap \Phi$ and $\mathring \Phi_J^+ = \bigoplus_{j \in J} \N \alpha_j \cap \mathring \Phi$, respectively.
\end{definition}
\begin{remark}
    Due to the way we defined $A$ from $\mathring A$, we have that the angles between $\abar_i, \abar_j$ and $\alpha_i,\alpha_j$ are the same. Thus, for every $J \in Q_A$ we have $\Phi_J \cong \mathring \Phi_J$.
\end{remark}
\begin{proposition} \label{prop: description image of phi}
We have the following description of the image of the local groups under $\phi$, where the product is taken in an arbitrary fixed order on $\Phi_J^+$.

\begin{align*}
    \phi(U_J) = \left\{ \prod_{a_{\alpha,m} \in \Phi_J^+} x_\alpha(\lambda_{a_{\alpha, m}} t^m)  \mid \lambda_{a_{\alpha,m}} \in k \right\}.
\end{align*}
\end{proposition}

\begin{proof}
    Note that, for $J \in Q_A $, $\mathcal{U}_J$ is the unitpotent subgroup of the positive standard Borel subgroup of $\chev_{(A_{i,j})_{i,j \in J}}(k)$. Thus, from \cite[Theorem 5.3.3]{carter1989simple} for height preserving orders and \cite[Lemma 16]{SteinbergChevalley} for arbitrary orders on $\Phi_J^+$, we get that 
    $$\mathcal{U}_J  = \left\{\prod_{\alpha \in \Phi_J^+} u_\alpha(\lambda_\alpha) \mid \lambda_\alpha \in k \right\}.$$
    Thus,
    \begin{align*}        
    \phi(U_J) &= \Tilde{\phi}(\mathcal{U}_J) = \left\{\prod_{\alpha \in \Phi_J^+} \Tilde{\phi}(u_\alpha(\lambda_\alpha)) \mid \lambda_\alpha \in k \right\} = \left\{\prod_{a_{\alpha,m} \in \Phi_J^+} \Tilde{\phi}(u_{a_{\alpha,m}}(\lambda_{a_{\alpha,m}})) \mid \lambda_{a_{\alpha,m}} \in k \right\} \\
    & = \left\{\prod_{a_{\alpha,m} \in \Phi_J^+} x_\alpha(\lambda_{a_{\alpha,m}}t^m) \mid \lambda_{a_{\alpha,m}} \in k \right\}.
    \end{align*}
\end{proof}

\begin{remark} \label{rmk: m=0 or 1}
    Let $a_{\alpha,m} \in \Phi_J^+$ for some $J \in Q_A$. If $0 \notin J$, then we can write $a_{\alpha,m} = \sum_{j \in J} n_j \abar_j = \overline{\sum_{j \in J} n_{j} \alpha_j}$ for some $n_j \in \N$ and thus $\alpha = \sum_{j \in J} n_{j} \alpha_j$ and $m = 0$. 
    
    If $0 \in J$, then, for some $n_j \in \N$, we have 
    $$a_{\alpha,m} = \sum_{j \in J} n_{j} \abar_j = \overline{\underbrace{\sum_{j \in J \setminus \{0\}} n_j \alpha_j}_{=\beta \in \mathring \Phi}} + n_0 \abar_0 = \overline{\beta} + \overline{-\gamma} + \overline{\delta} = \overline{ \beta - \gamma} + \overline{\delta}$$
    and thus $m = 1$.
\end{remark}

\begin{lemma} \label{lemma: inj on lo grps affine case}
    Let $f \in k[t]$ be an irreducible polynomial with $\deg(f) \geq 2$. Then the map $\varphi_f = \pi_f \circ \phi: \uak \to \chev_{\Ao}(k[t]/(f))$ is injective on the local groups $U_J, J \in Q_A$.
\end{lemma}
\begin{proof}
    We know that $\phi$ is injective on the local groups, since $\psi$ is injective on the local groups and $\Tilde{\phi}$ is injective. Since $\deg(f) \geq 2$ we have that $\pi_f|_k$ and $\pi_f|_{kt}$ are injective. The results of \Cref{prop: description image of phi} and \Cref{rmk: m=0 or 1} yield that the entries in image of a local group are polynomials of degree less or equal to 1, thus $\pi_f$ is injective on the image of the local groups. 
\end{proof}
Next, we show that $\varphi_f$ satisfies the intersection property (\ref{essential property}) of \Cref{section: main theorem}.
\begin{lemma}
    Let $k$ be a finite field. Then there exists $\ell \in \N$, depending only on the GCM $A$ and $k$, such that for any irreducible polynomial $f \in k[t]$ with $\deg(f) \geq \ell$ and for all $J,K\in Q_A$, we have
    $$\varphi_f(U_J) \cap \varphi_f(U_K) = \varphi_f(U_J \cap U_K).$$
\end{lemma}
\begin{proof}
We know that in $\phi(U_J), \phi(U_K)$ the degree of the entries inside the $x_\alpha$ is at most one by \Cref{prop: description image of phi}. Since the sets $U_J, U_K$ are finite, we know that the degree of $t$ inside the $x_\alpha$ that appear in $\phi(U_J) \cdot \phi(U_K)$ is bounded above by some constant $\ell_{J,K}$ depending on the size and the type of $J,K$. Set $\ell = \max_{J,K \in Q_A} \ell_{J,K} +1$. Then, by choosing $f$ of degree at least $\ell $, we ensure that $\pi_f$ is injective on $\phi(U_J) \cdot \phi(U_K)$ for all $J,K \in Q_A$.

Furthermore, we know that $\Tilde{\phi}$ is injective. Thus $\pi_f \circ \Tilde{\phi}$ is injective on $\bigcup_{J,K \in Q_A} \psi(U_J) \cdot \psi(U_K)$. Hence, \Cref{lemma: composition and stupid property} gives the desired result.
\end{proof}

For fields of characteristic different from 2, this result can be made substantially more precise using the adjoint representation of $\chev_{\Ao}(k[t]/(f))$. For the characteristic 2 case, see also \Cref{rmk: characteristic 2}.
\begin{proposition} \label{prop: affine case IP for characteristic not 2 and 3}
    Let $k$ be a finite field of characteristic different from 2 and let $f \in k[t]$ be irreducible with $\deg(f) \geq 2$. Then for all $J,K\in Q_A$, we have
    $$\varphi_f(U_J) \cap \varphi_f(U_K) = \varphi_f(U_J \cap U_K).$$
\end{proposition}

\begin{proof}
Recall that 
\begin{alignat*}{2}
    \varphi_f: \uak & \overset{\phi}{\longrightarrow} \chev_{\Ao}(k[t]) && \overset{\pi_f}{\longrightarrow} \chev_{\Ao}(k[t]/(f))\\
    u_i(\lambda) & \mapsto \begin{cases} x_{\alpha_i}(\lambda)  & i \neq 0 \\ x_{-\gamma}(\lambda t) & i = 0 \end{cases} && \mapsto \begin{cases} x_{\alpha_i}(\lambda + (f))  & i \neq 0 \\ x_{-\gamma}(\lambda t + (f)) & i = 0 \end{cases}.
\end{alignat*}

We show that the map $\pi_f$ is injective on $\phi(U_J) \cdot \phi(U_K)$. To do so, we consider the adjoint representation of the Chevalley group $\Ad: \chev_{\Ao}(k[t]/(f)) \to \operatorname{GL}(\mathfrak{g})$ and show that the composition $\Ad \circ \pi_f$ is injective on $\phi(U_J) \cdot \phi(U_K)$. Here $\mathfrak{g}$ denotes the Lie algebra of $G = \chev_{\Ao}(k[t]/(f))$.

Let $\mathfrak{h}$ denote a Cartan subalgebra of $\mathfrak{g}$, and let $\mathfrak{g}_\alpha = \{x \in \mathfrak{g} \mid [h,x] = \alpha(h) x \text{ for all } h \in \mathfrak{h} \}$ be the eigenspace corresponding to $\alpha \in \mathfrak{h}^*$. We can identify the elements $\alpha \in \mathfrak{h}^*$ with $\mathfrak{g}_\alpha \neq 0$ with the elements of $\mathring \Phi$. This gives rise to a root space decomposition of $\mathfrak{g}$:
$$\mathfrak{g} = \mathfrak{h} \oplus \bigoplus_{\alpha \in \mathring \Phi} \mathfrak{g}_\alpha.$$

Each $\mathfrak{g}_\alpha$ is one-dimensional and we can choose a basis vector $e_\alpha$ of $\mathfrak{g}_\alpha$ such that $\exp(se_\alpha) = x_\alpha(s)$ for every $s \in k[t]/(f) =:K_f$. Set $h_\alpha = [e_\alpha, e_{-\alpha}] \in \mathfrak{h}$ for $\alpha \in \mathring \Phi$. 

In $\mathfrak{g}$ the following relations hold, see \cite[Chapter 4.2]{carter1989simple}:
\begin{align*}
    [h_\alpha,h_\beta] &= 0 \text{ for all } \alpha,\beta \in \mathring \Pi; \\
    [h_\alpha,e_\beta] &= C_{\alpha,\beta}e_\beta \text{ for all } \alpha \in \mathring \Pi, \beta \in \mathring \Phi; \\
    [e_\alpha,e_\beta] &= \begin{cases}
    0 & \text{ if } \alpha+\beta \notin \mathring \Phi \\ N_{\alpha,\beta} e_{\alpha+\beta} & \text{ if } \alpha+\beta \in \mathring \Phi \end{cases} \text{ for all } \alpha,\beta \in \mathring \Phi, \alpha \neq -\beta ;\\
    [e_\alpha, e_{-\alpha}]  &= h_\alpha \text{ for all } \alpha \in \mathring \Phi.
\end{align*}
Here $C_{\alpha,\beta}$ and $N_{\alpha,\beta}$ are integer constants given in the following way: $C_{\alpha,\alpha} = 2$, and, for $\alpha \neq -\beta$, let $p,q\in \N$ be such that for all $-p \leq i \leq q$ we have $i \alpha + \beta \in \mathring \Phi$ but $-(p+1) \alpha + \beta \notin \mathring \Phi, (q+1)\alpha + \beta \notin \mathring \Phi$. Then $C_{\alpha,\beta} = p-q$ and $N_{\alpha, \beta} = \pm (p+1)$ where the sign depends on the explicit choice of $e_\alpha, e_\beta$. Since $\mathring \Phi$ is spherical we have that $p,q \leq 3$ and hence $|C_{\alpha,\beta}| \leq 3, |N_{\alpha, \beta}| \leq 4$.

The adjoint representation $\Ad: G \to \operatorname{GL}(\mathfrak{g})$ has the following concrete realization \cite[3.2.5]{BruhatTits84}:
\begin{alignat}{2}
    (\Ad x_\alpha(s))(e_\beta) &= \sum_{j \geq 0, \beta + j \alpha \in \mathring \Phi} M_{\alpha,\beta,j} s^je_{\beta + j \alpha} \qquad && \text{ for } \alpha, \beta \in \mathring \Phi, \alpha \neq - \beta, s \in K_f \label{ad rep: different ealpha}\\
    (\Ad x_\alpha(s))(e_{-\alpha}) &=e_{-\alpha} - s h_\alpha + s^2 e_\alpha && \text{ for } \alpha \in \mathring \Phi, s \in K_f\\
    (\Ad x_\alpha(s))(h) &=h- \alpha(h) s e_\alpha && \text{ for } \alpha \in \mathring \Phi, h \in \mathfrak{h}, s \in K_f. \label{ad rep: x alpha and h}
\end{alignat} 
Here the $M_{\alpha, \beta,j}$ are constants given by the structure constants (from the commutator relations in the Chevalley group) $C^{\alpha,\beta}_{i,j}$ in the following way: $M_{\alpha,\beta, j} = C^{\alpha, \beta}_{j, 1}$ for $j \geq 1$ and $M_{\alpha,\beta, 0} = 1$.

By \cite[Theorem 7.55]{marquis2018introduction} there exists a maximal split torus $T$ of $G$ such that $\ker(\Ad) \subseteq T$. Furthermore, we have $\pi_f(\phi(U_J)) \cap T = \pi_f(\phi(U_K) )\cap T = \{1\}$. Hence, $\Ad$ is injective on $\pi_f(\phi(U_J))$ and $\pi_f(\phi(U_K))$.

Next, we show that $\ker(\Ad) \cap \pi_f(\phi(U_J) \phi(U_K)) = \{1\}$.
Suppose there exists a $1\neq z \in \pi_f(\phi(U_J)\phi(U_K))$ that acts trivially on $\mathfrak{g}$, i.e. $\Ad(z) = 1$. Then there exist $g \in \pi_f(\phi(U_J)), g' \in \pi_f(\phi(U_K))$ such that $\Ad(g) = \Ad(g')$. Without loss of generality, we can assume that $J = I \setminus \{i\}$ for some $i\in I$ (otherwise $J \subset I\setminus \{i\}$ for some $i\in I$, thus $U_J \subseteq U_{I\setminus \{i\}}$ and hence $\pi_f(\phi(U_J)) \subseteq \pi_f(\phi(U_{I\setminus\{i\}}))$). Assuming that $K \nsubseteq J$ (otherwise the statement is trivial) we get that $i \in K$.

Let $L_J$ be the subspace of $\mathfrak{g}$ spanned by the basis $\{h_\alpha,  e_\beta \mid \alpha\in \mathring \Pi, \beta \in \mathring \Phi^+_J \}$. Then $L_J$ is a Lie subalgebra of $\mathfrak{g}$. Note that $\Ad(g)$ stabilizes $L_J$ and that, if we define $L_K$ analogously to $L_J$, we have that $\Ad(g')$ stabilizes $L_K$.

Next, we show that there exists a $h \in \mathfrak{h}$ with $\Ad(g')(h) \notin L_J$.

Note that any $\alpha \in \ophiK$ has a unique decomposition as $\sum_{i \in K} n_i \alpha_i$ for some $n_i \in \N$. We define the \textit{height} of a root $\alpha$ in $\ophiK$ to be $\height(\alpha) = \sum_{i \in K} n_i$.
We order the roots in $\ophiK$ in the following way: we start with the roots in $\ophiK \setminus \ophiJ$ and order them increasing in height. Afterwards, we arrange the remaining roots, which are in $\ophiK \cap \ophiJ$, in arbitrary order. 
Then $g'$ can be written in the following way, where the product is taken in the order that we just fixed (see \Cref{prop: description image of phi}):
$$g' = \prod_{\alpha \in \mathring \Phi_K^+}x_\alpha(s_\alpha).$$
Here the $s_\alpha$ are elements of $K_f$. Let $\alpha \in \ophiK$ be the first root in the order such that $s_\alpha \neq 0$. Since $g' \notin \pi_f(\phi(U_J))$, we have that $\alpha \in \ophiK \setminus \ophiJ$.

The set $P:=\ophiK \setminus \ophiJ$ has the following two properties. For all $\beta, \gamma \in P$ such that $\beta+\gamma \in \mathring \Phi ^+$ we have $\beta+\gamma \in P$ (we say that $P$ is \textit{closed}), and for all $\beta \in P , \gamma \in \ophiK$ with $\beta + \gamma \in \mathring \Phi^+ $ we have $\beta + \gamma \in P$ (we say that $P$ is an ideal in $\ophiK$).

%??? Note that, for $\beta, \gamma \in \ophiK$ with $\beta + \gamma \in \ophiK$, we have $\height(\beta + \gamma) = \height(\beta)+\height(\gamma)$ thus $\alpha$ cannot be written as a sum of roots from $\{\beta \in \ophiK \mid s_\beta \neq 0\}$.

We set $h= h_\alpha$, the co-root corresponding to $\alpha$.
We inductively apply formulas (\ref{ad rep: different ealpha}) and (\ref{ad rep: x alpha and h}) to $\Ad(g')(h_\alpha) = \Ad(\prod_{\alpha \in \mathring \Phi_K^+}x_\alpha(s_\alpha))(h_\alpha)$, using the observations above and that $x_\alpha(s_\alpha)$ is the first non-trivial term in $g'$. Further using that $\Ad(x_\alpha(s_\alpha))(h_\alpha) = h_\alpha -2s_\alpha e_\alpha$, we have that 
$$\Ad(g')(h_\alpha) = h_\alpha -2s_\alpha e_\alpha + \sum_{\beta \in \ophiJ \cap \ophiK} \lambda_\beta e_\beta + \sum_{\gamma \in \ophiK \setminus \ophiJ: \ \height(\gamma) \geq \height(\alpha), \gamma \neq \alpha} \lambda_\gamma e_\gamma$$
for some $\lambda_\beta, \lambda_\gamma \in K_f$.

Since we assume $\operatorname{char}(k)\neq 2$, we have $\Ad(g')(h_\alpha) \notin L_J$, a contradiction to $\Ad(g) = \Ad(g')$. Therefore, we have shown that $\ker(\Ad) \cap \pi_f(\phi(U_J) \phi(U_K)) = \{1\}$.

Now, we can conclude that $\pi_f$ is injective on $\phi(U_J)\phi(U_K)$. Let $1 \neq a\in \phi(U_J), 1 \neq b \in \phi(U_K)$ such that $ab \neq 1$. We want to show that $\pi_f (ab) \neq 1$. By \Cref{lemma: inj on lo grps affine case}, we have $\pi_f(a) \neq 1, \pi_f(b) \neq 1$. Since both are in $\pi_f(\phi(U_J) \phi(U_K))$ and $\pi_f(\phi(U_J) \phi(U_K)) \cap \ker (\Ad) = \{1\}$, we know $\Ad(\pi_f(a)) \neq 1, \Ad(\pi_f(b))\neq 1$ and $\Ad(\pi_f(a)\pi_f(b))= \Ad(\pi_f(ab)) \neq 1$. This implies that $\pi_f(ab) \neq 1$. 

Furthermore, we know that $\Tilde{\phi}: U^+ \to \chev_{\Ao}(k[t])$ is injective. Thus $\pi_f \circ \Tilde{\phi}$ is injective on $\bigcup_{J,K \in Q_A} \psi(U_J) \cdot \psi(U_K)$. Hence, \Cref{lemma: composition and stupid property} gives the desired result.
\end{proof}

\begin{remark} \label{rmk: characteristic 2}
    If the characteristic of $k$ is 2, the proof of \Cref{prop: affine case IP for characteristic not 2 and 3} only fails if $\alpha(h)=0$ for all $h$ in the Cartan subalgebra. The root $\alpha$ is the image of a simple root, say $\alpha_i$, under the action of an element of the Weyl group. Thus, 
    $$\alpha(h) = 0 \text{ for all } h \in \mathfrak{h} \iff \alpha_i(h) = 0 \text{ for all } h \in \mathfrak{h} \iff \mathring A_{j,i} = \alpha_i(h_j) = 0 \text{ for all } j \in \mathring I$$
    If the Cartan matrix $\mathring A$ has no column equal to 0 mod 2, i.e. the Dynkin diagram has no double edges, then the above proof still works by choosing $h$ such that $\alpha(h) \neq 0$ in $k$ instead of taking $h_\alpha$. Hence, if the characteristic of $k$ is 2 and the Cartan matrix is not of type $B_n$ or $C_n$, we still get the result of \Cref{prop: affine case IP for characteristic not 2 and 3}.
\end{remark}

Finally, we show that $\varphi_f$ is surjective.

\begin{proposition} \label{prop: image of phi}
    The image of $\phi$ has the following form:
    $$\phi(\uak) = \left\langle x_\alpha(g) \mid \alpha \in \mathring \Phi, g \in k[t] \text{ if } \alpha \in \mathring \Phi^+, g \in (t) \lhd k[t] \text{ otherwise } \right\rangle.$$
\end{proposition}
\begin{proof}
    Recall that $\Phi^+ = \{\abar + m\overline{\delta} \mid m = 0 \text{ and } \alpha \in \mathring \Phi^+ \text{ or } m \geq 1, \alpha \in \mathring \Phi \}$ and for $u_{a_{\alpha,m}}(\lambda) \in U^+$ we have $\Tilde{\phi}(u_{a_{\alpha,m}}(\lambda)) = x_\alpha(\lambda t^m)$. Thus, 
\begin{align*}
\phi(\uak) &= \Tilde{\phi}(U^+) = \Tilde{\phi}(\langle u_{a_{\alpha,m}}(\lambda_{a_{\alpha,m}}) \mid a_{\alpha,m} \in \Phi^+, \lambda \in k \rangle) \\
&= \langle x_\alpha(\lambda_\alpha t^{m_\alpha}) \mid m = 0 \text{ and } \alpha \in \mathring \Phi^+ \text{ or } m \geq 1, \alpha \in \mathring \Phi  \rangle \end{align*}
    which implies the result.
\end{proof}

\begin{corollary}
     The map $\varphi_f = \pi_f \circ \phi: \uak \to \chev_{\Ao}(k[t]/(f))$ is surjective. 
\end{corollary}
\begin{proof}
    This follows from \Cref{prop: image of phi} together with the fact that $\pi_f((t)) = \pi_f(k[t])$.
\end{proof}

The work done in this section results in the following theorem.

\begin{theorem}
    Let $\mathring A$ be an irreducible Cartan matrix and let $k$ be a finite field. Then there exists $\ell \in \N$, where $\ell = 2$ if $k$ has characteristic different from 2, such that the following holds. Any family of irreducible polynomials  $f_m \in k[t], m \in \N$, satisfying $\deg(f_m) \geq \ell$ for all $m \in \N$ and $\deg(f_m) \to \infty$ for $m\to \infty$, gives rise to an infinite family of maps $(\varphi_m)_{m \in \N}$ satisfying all requirements of \Cref{thm: main thm} and such that $|\operatorname{Im}(\varphi_m)| = |\chev_{\Ao}(\mathbb{F}_{p^{\deg(f_m)}})| \overset{m \to \infty}{\longrightarrow} \infty$. Thus, each such family of irreducible polynomials gives rise to an infinite family of bounded degree HDX.
\end{theorem}

\subsection{Example case: $\operatorname{SL}_3$}\label{subsection: quot in Chev: sl3}
In this section, we describe the situation in the case where $\mathring \Phi$ is of type $A_2$, i.e. $\chev_{\Ao} = \operatorname{SL}_3$.

This will give rise to probably the most basic setting to which we can apply our main result \Cref{thm: main thm}.

In this case the generators $x_\alpha(\lambda)$ for $\alpha \in \mathring \Phi$ of the Chevalley group can be realized as the following matrices:
\begin{align*}
    x_{\alpha_1}(\lambda) = \begin{pmatrix} 1 & \lambda & 0 \\ 0 & 1 & 0 \\ 0 & 0 & 1\end{pmatrix}, \qquad x_{\alpha_2}(\lambda) = \begin{pmatrix} 1 & 0 & 0 \\ 0 & 1 & \lambda \\ 0 & 0 & 1\end{pmatrix}, \qquad x_{-\gamma}(\lambda) = \begin{pmatrix} 1 & 0 & 0 \\ 0 & 1 & 0 \\ \lambda & 0 & 1\end{pmatrix}.
\end{align*}

We consider the map $\phi$ as described above:
\begin{align*}
    \phi: \uak = \langle u_i(\lambda) \mid i=0,1,2, \lambda \in k \rangle &\to \chev_{\Ao}(k[t]) = \langle x_\alpha(g) \mid \alpha \in \mathring \Phi, g \in k[t] \rangle \\
    u_1(\lambda)  \mapsto x_{\alpha_1}(\lambda) \qquad     u_2(\lambda) &\mapsto x_{\alpha_2}(\lambda) \qquad    u_0(\lambda) \mapsto x_{-\gamma}(\lambda t).
\end{align*}

We get the following images of the local groups of the spherical subsets of $I = \{0,1,2\}$:
$$\phi(U_{12})= \begin{pmatrix} 1 & k & k \\ 0 & 1 & k \\ 0 & 0 & 1\end{pmatrix}, \qquad \phi(U_{01})= \begin{pmatrix} 1 & k & 0 \\ 0 & 1 & 0 \\ kt & kt & 1\end{pmatrix}, \qquad \phi(U_{02})= \begin{pmatrix} 1 & 0 & 0 \\ kt & 1 & k \\ kt & 0 & 1\end{pmatrix}.$$

Next, we fix a family of irreducible polynomials $(f_m)_{m \in \N}$ with $\deg(f_m) \geq 2$ for all $m \in  \N$, and with $\lim_{m \to \infty} \deg(f_m) = \infty$. Note that when $\mathring \Phi$ is of type $A_n$, then the constants $C_{\alpha,\beta}$ appearing in the proof of \Cref{prop: affine case IP for characteristic not 2 and 3} have absolute value at most one and thus do not vanish when seen in $k[t]/(f_m)$ even if $k$ has characteristic 2. This gives rise to maps
\begin{align*}
    \varphi_m = \pi_{f_m} \circ \phi: \uak \to \chev_{\Ao}(k[t]/(f_m)) 
\end{align*}
by taking each entry of the matrices modulo the ideal generated by $f_m$.

Setting 
$$G_m = \chev_{\Ao}(k[t]/(f_m)), \qquad H_0^m = \pi_{f_m}(\phi(U_{12})), \qquad H_1^m = \pi_{f_m}(\phi(U_{02})), \qquad H_2^m = \pi_{f_m}(\phi(U_{01}))$$ 
we get an infinite family of HDX of the form $X_m = \CC(G_m; H_0^m, H_1^m, H_2^m), m \in \N$.

\section{Comparison and generalisation}\label{section : context}
In this section, we first compare our construction to those of Kaufman and Oppenheim in \cite{kaufman2018construction} and O'Donnell and Pratt in \cite{o2022high}. Then, we give a potentially more general formulation of \Cref{thm: main thm} that might give rise to even more high-dimensional expanders in the future. 

\subsection{Comparison to other constructions}\label{subsection: quot in chev: comparison}
The first construction of HDX using coset complexes by Kaufman and Oppenheim \cite{kaufman2018construction} takes as groups the elementary groups $G_{m} = El_{n+1}(\mathbb{F}_q/(t^m))$ and subgroups $H_i = \langle e_{j,j+1} (a+bt) \mid j \in \{0,\dots,n\} \setminus \{i\}, a,b \in \mathbb{F}_q \rangle$ where $j+1$ is considered modulo $n+1$ and $e_{j,j+1}(g)$ is the matrix with ones on the diagonal and $g$ at position $(j,j+1)$.

On the other hand, O'Donnell and Pratt take in \cite{o2022high} as groups $G_m = \chev_{\Ao}(\mathbb{F}_{p^m})$ and as subgroups $H_\alpha = \langle x_\beta(g) \mid \beta \in S\setminus \{\alpha\}, g\in \mathbb{F}_{p^m}, \deg(g) \leq 1 \}$ where $S= \mathring \Pi \cup \{ -\sum_{i=1}^d \alpha_i\}$ and $\alpha \in S$. \\

The construction by Kaufman and Oppenheim can be seen as an instance of O'Donnell and Pratt's construction.

These constructions differ from our construction applied to Chevalley groups in two places: they allow the entries of the generators to be polynomials of degree $\leq 1$ while in our construction the entries will either be constant or multiples of $t$ (in both cases modulo an irreducible polynomial). Secondly, they take as extra generator the one indexed by the root $-\sum_{i=1}^d \alpha_i$ while we take the negative of the highest root. In case the root system is of type $A_n$ (and the Chevalley group is thus $\operatorname{SL}_{n+1}$) those two definitions coincide, but not in the other cases. The authors note in \cite[Remark 3.11]{o2022high} that other choices of $S$ would be possible. Finally, another difference is that our construction works in the $G_2$ case, and \Cref{thm: main thm} gives a tool that allows to consider HDX over many other groups as well.

The coset complexes described in \Cref{subsection: quot in Chev: sl3} are very similar to the coset complexes used in \cite{DinurNewCodesonHDX} to construct a family of symmetric error-correcting codes with low-density parity-check matrices. But attempts of showing that the appearing complexes are isomorphic have failed so far. 

Thus, our construction is not a precise generalization of the previous constructions, but it seems that the underlying structure that leads to the good expansion property is the same. The way in which the bound of the degree is achieved differs. Our main result, \Cref{thm: main thm}, utilizes the underlying structure that leads to good expansion, in the most general way so far.

\subsection{Generalization of the main theorem}\label{subsection : generalization}

We can formulate \Cref{thm: main thm} in a potentially more general context than KMS-groups, although it is not clear whether there are other groups that satisfy the requirements of the theorem and give rise to interesting HDX.

\begin{theorem} \label{thm: generalized main thm}
   Let $I$ be a finite set and for each $J \subset I, 1 \leq |J| \leq 2$ let finite groups $U_J$ be given together with inclusions $U_{\{i\}} \hookrightarrow U_{\{i,j\}}$ for all $i,j \in I$. Furthermore, we require that $U_{\{i,j\}} = \langle U_{\{i\}}, U_{\{j\}}\rangle $ for all $i,j \in I$. Let $U$ be the direct limit of the $U_J, J \subset I, 1 \leq |J| \leq 2$ with respect to the inclusions, i.e. $U = *_{J\subset I, 1 \leq |J| \leq 2} U_J/(U_{\{i\}} \hookrightarrow U_{\{i,j\}})$. We require that 
   \begin{enumerate}[label=(\roman*)]
    \item there exists a $0<\lambda \leq \frac{1}{|I|-1}$ such that $\CC(U_{\{i,j\}};(U_{\{i\}},U_{\{j\}}))$ is a $\lambda$-expander graph for all $i \neq j \in I$;
   \item the groups $U_{\{i,j\}}$ must embed into $U$, i.e. the natural map from $U_{\{i,j\}}$ to $U$ must be injective. Then we identify $U_{\{i\}}, U_{\{i,j\}}$ with their images in $U$.
\end{enumerate}
Set $U_J:= \langle U_{\{j\}} \mid j \in J \rangle \leq U$ for any $J \subsetneq I$. We further require that
\begin{enumerate}     
       \item[(iii)] $|U_J| < \infty$ for all $J \subsetneq I$;
       \item[(iv)] $U_J \cap U_K = U_{J\cap K}$ for all $J,K \subsetneq I$.
   \end{enumerate}
    Let $G$ be a finite group such that there exists a surjective homomorphism $\phi: U \to G$ that is injective on all the $U_J$ and such that
    $$\phi(U_J) \cap \phi(U_K) = \phi(U_J \cap U_K) \text{ for all } J,K \subsetneq I.$$
Set $H_i = \phi(U_{I\setminus \{i\}})$, then $\CC(G; (H_i)_{i \in I})$ is a $\gamma$-spectral high-dimensional expander, for some $\gamma>0$ independent of $G$ and $\phi$. An infinite family of such groups $G$ and maps $\phi$, where the size of the groups tends to infinity, gives rise to an finite family of bounded degree HDX.
\end{theorem}

Proving this theorem works step by step exactly like the proof of \Cref{thm: main thm}, since we assumed exactly those properties of KMS-groups that we needed. One can interpret \Cref{section: KMS group and root subgroups} as a section in which we prove that KMS groups satisfy the assumptions of \Cref{thm: generalized main thm}. 

For $|I| = 3$, there exist other groups $U_{\{i\}}, U_{\{i,j\}}$ that satisfy the conditions of \Cref{thm: generalized main thm}, for example the ones constructed in \cite[Theorem 4.6]{LubGeneralizedTriangleGrps} or \cite[Chapter 4, Example 1]{RonanLecturesonBuildings}. In both examples, the difficulty is to find appropriate families of finite quotients which would be necessary to obtain an infinite family of bounded degree HDX.

%----------------------
%%%% Bibliography
{
\renewcommand*\MakeUppercase[1]{#1}
\bibliographystyle{alpha}							
\bibliography{references}						
}

%---------------

\end{document}